\documentclass[11pt,reqno,a4paper]{amsart}
\usepackage{graphicx}
\usepackage{amsmath,amsfonts,amssymb,amsrefs,cases}
\usepackage{xcolor}
\usepackage{dsfont}
\usepackage{etoolbox}
\usepackage{esint}

\usepackage[tmargin=2.7cm,bmargin=2.7cm,hmargin=2.6cm,bindingoffset=-0cm]{geometry}
\usepackage{titlesec}
\usepackage{titletoc}
\usepackage{pgf,tikz}

\usepackage[showboxes,overlay]{textpos}

\TPshowboxesfalse
\usepackage{graphicx}

 \usepackage[colorlinks=true,linkcolor=black,citecolor=blue]{hyperref}
\usepackage{ulem}

\usepackage[utopia]{mathdesign}

\usetikzlibrary{arrows}

\DeclareMathOperator{\rt}{ curl}

\DeclareMathOperator{\dv}{div}
\DeclareMathOperator{\I}{I_3}

\renewcommand{\Subset}{\subsubset}

\newcommand{\mP}{\mathcal{P}}

\newcommand{\Y}{ \Gamma  }

\newcommand{\muef}{\bmu^{\mathrm{eff}}}
\newcommand{\epsef}{\beps^{\mathrm{eff}}}

\newcommand{\R}{\mathbb{R}}

\newcommand{\N}{\mathbb{N}}
\newcommand{\Z}{\mathbb{Z}}
\newcommand{\C}{\mathbb{C}}

\newcommand{\un}{\mathds{1}}

\newcommand{\Hunper}{W^{1,2}_\sharp}
\newcommand{\Hun}{W^{1,2}}

\newcommand{\cvd}{\rightharpoonup\hspace{-0.3cm}\rightharpoonup}
\newcommand{\cvdf}{\rightarrow\hspace{-0.3cm}\rightarrow}
\newcommand{\cvw}{\rightharpoonup}

\newcommand{\rp}{\right)}
\newcommand{\lp}{\left(}
\newcommand{\ov}{\overline}
\newcommand{\rc}{\right]}
\newcommand{\lc}{\left[}

\newcommand{\ra}{\right\}}
\newcommand{\la}{\left\{}
\renewcommand{\ln}{\left\|}
\newcommand{\rn}{\right\|}

\newcommand{\x}{\times }
\newcommand{\dom}{\Omega}
\newcommand{\im}{\Im}

\newcommand{\subsubset}{\subset\!\subset}
\newcommand{\eps}{\varepsilon}
\newcommand{\disp}{\displaystyle}
\newcommand{\xxeta}{\Big( x,\frac{x}{\eta}\Big)}
\newcommand{\xeta}{\Big( \frac{x}{\eta}\Big)}
\newcommand{\med}{\medskip\noindent}

\newcommand{\bfu}{\boldsymbol{u}}

\newcommand{\bfv}{{\boldsymbol v}}

\newcommand{\bfc}{{\boldsymbol c}}

\newcommand{\bff}{{\boldsymbol f}}
\newcommand{\bfg}{{\boldsymbol g}}
\newcommand{\bfn}{{\boldsymbol n}}
\newcommand{\bfw}{{\boldsymbol w}}

\newcommand{\bfy}{{\boldsymbol y}}
\newcommand{\bfz}{{\boldsymbol z}}
\newcommand{\bfe}{{\boldsymbol e}}
\newcommand{\bfk}{{\boldsymbol k}}

\newcommand{\Einc}{{\boldsymbol E}^{\mathrm{inc}}}
\newcommand{\Hinc}{{\boldsymbol H}^{\mathrm{inc}}}

\newcommand{\bmu}{{\boldsymbol \mu}}
\newcommand{\beps}{{\boldsymbol \eps}}

\newcommand{\bs}{ \boldsymbol }
\newcommand{\bse}{ \boldsymbol{E} }
\newcommand{\bsh}{ \boldsymbol{H} }

\newcommand{\bseinc}{ \boldsymbol{E}^\mathrm{inc} }
\newcommand{\bshinc}{ \boldsymbol{H}^\mathrm{inc} }

\newcommand{\bsB}{ \boldsymbol{B} }

\newcommand{\bsb}{ \boldsymbol{b} }

\newcommand{\bsj}{ \boldsymbol{J} }
\newcommand{\bsu}{ \boldsymbol{u} }
\newcommand{\bsv}{ \boldsymbol{v} }

\newcommand{\bsx}{ \boldsymbol{x} }
\newcommand{\bpsi}{ \boldsymbol{\psi} }
\newcommand{\bphi}{ \boldsymbol{\varphi} }
\newcommand{\nshape}[1]{ {\upshape #1}}

\newcommand{\dx}{ \,\mathrm{d}x }
\newcommand{\dy}{ \,\mathrm{d}y }

\newtheorem{theorem}{Theorem}[section]
\newtheorem{proposition}[theorem]{Proposition}
\newtheorem{lemma}[theorem]{Lemma}

\newtheorem{definition}[theorem]{Definition}
\newtheorem{remark}[theorem]{Remark}
\numberwithin{equation}{section}

\titleformat{\section}[block]{\bfseries\large\filcenter\scshape}{\thesection.}{0.4em}{}
\titleformat{\paragraph}[runin]{\bfseries}{}{}{}
\titlespacing*{\paragraph}{0pt}{0.5\baselineskip}{0.5\baselineskip}

\titleformat{\subsection}[block]{\bfseries}{\thesubsection.}{0.4em}{}
\titlespacing*{\subsection}{0pt}{\baselineskip}{0.1\baselineskip}

\dottedcontents{section}[3.5em]{\bfseries\scshape}{2em}{1pc}

\dottedcontents{subsection}[6.1em]{\small}{2.5em}{1pc}

\author{Guy Bouchitt\'e, Christophe Bourel, Didier Felbacq}
\address[G. Bouchitt\'e]{IMATH, Universit\'e du Sud-Toulon-Var, 83957 La Garde cedex, France}
\address[C. Bourel]{Univ. Littoral Côte d'Opale, EA 2797 - LMPA, F- 62228 Calais, France}
\address[D. Felbacq]{L2C, Universit\'e de Montpellier,  34000 Montpellier, France}

\begin{document}

\title{Homogenization near resonances and artificial magnetism in 3D dielectric metamaterials} 
 

\begin{abstract}

It is now well established that the homogenization of a periodic array of parallel dielectric fibers with suitably scaled high permittivity
can lead to a (possibly) negative frequency-dependent effective permeability. However this result based on a two-dimensional approach holds merely in the case of linearly polarized magnetic fields, reducing thus its applications to
infinite cylindrical obstacles.
 In this paper we consider a dielectric structure placed in a \textit{bounded} domain of $\R^3$ 
and perform a full 3D asymptotic analysis.
The main ingredient is a new averaging method for characterizing the bulk effective magnetic field in the vanishing-period limit.
 We evidence a vectorial spectral problem on the periodic cell which determines micro-resonances and 
encodes the oscillating behavior of the magnetic field from  which  artificial magnetism arises.
At a macroscopic level we deduce an effective permeability  tensor that we can be make explicit as a function of the frequency. 
As far as sign-changing permeability are sought after, we may foresee that periodic bulk dielectric inclusions could be an
efficient alternative to the very popular metallic split-ring structure proposed by Pendry. 
Part of these results have been announced in~\cite{bou_bou_cras}.

\end{abstract}

\maketitle


\tableofcontents

\section{Introduction and description of the model}

\bigskip

\paragraph{Physical background and recent mathematical progress.}

The behavior of an homogeneous material with respect to electromagnetic waves is characterized by its electric permittivity $\eps(\omega)$ and its magnetic permeability $\mu(\omega)$ : two physical quantities which depend on the frequency $\omega$.
For frequencies of visible light, the permittivity is a complex number $\eps:=\eps'+i\eps''$ where $\eps'\in\R$ (can be negative for some metals) and $\eps''\ge0$ (for passive media and if the harmonic
time-dependence is assumed to be $e^{-i\omega t}$). It is different for the permeability because all natural materials present a non-magnetic behavior in the visible region of the spectrum, i.e. their relative permeability is very close to one as in vacuum.

For the past fifteen years, there have been many researches on the realization of artificial materials, generally periodically micro-structured, behaving as homogeneous media, i.e. described by effective tensors $\epsef$ and $\muef$. An important issue is to design structures which allows a non-trivial permeability (possibly negative), a negative permittivity or both. The later case corresponds to a ``left-handed medium'' presenting a \textit{negative refractive index}. 

\bigskip
The first metamaterial possessing a negative effective permittivity was proposed by Pendry in 1996  \cite{pendry} and consists
in high conductivity parallel fibers occupying a very small volume fraction. A rigorous proof of this effective behavior, based on  homogenization techniques, appeared in \cites{pendry,boufel_random_media,Bloch_vector} in the case of infinitely long fibers and under a polarization assumption.
Surprisingly the same kind of behavior does not hold for finite-length fibers as demonstrated in\cites{boufel_cap,josa} where the
resulting permittivity law  is shown to be non-local. However  by inserting such a finite structure in  a larger scale structure,
 a reiterated homogenization procedure makes it possible mathematically to reach effective tensors $\epsef$ with negative eigenvalues (see \cite{cicp}).

\bigskip
In a similar way it is a challenging issue to design metamaterials able to display an artificial magnetic activity i.e. $\muef(\omega) \not=1$
and more specifically such that $\Re(\muef(\omega))<0$ in some range of frequencies. In photonic devices such a property
is usally explained by the ability of the structure to induce a local magnetic field presenting a Fano-like resonance \cite{Fano}. 
 The first and the most famous structure illustrating this phenomenon was proposed  by Pendry in \cites{pendry,pendry2} and consists in a periodic set of \textit{metallic split-ring resonators}. 
A field incident on this device induces micro-currents looping in each ring from which results  a macroscopic magnetic moment.
The mathematical study of this structure was made recently in \cite{kohn_ship} for the 2D case and in \cite{bou_schw} for the 3D general case. 

\medskip
Another way to produce artificial magnetism from dielectric structures was proposed in \cite{Pendry_mu}. Therein another kind of internal resonances is exploited, the so-called \textit{ Mie resonances}. These resonances take place inside each  dielectric inclusion and  generate  loops of displacement current inside the obstacle.
 This phenomenon can be evidenced experimentally on a composite structure  with a much simpler geometry than the one of the split-rings: it consists of
 periodically disposed  micro-cavities filled with a high dielectric material \cite{exp_mu_neg}.
 Subsequent works have shown the interest of Mie resonances to tailor the properties of dielectric metamaterials. These allow for the control of the Purcell effect \cite{xiao2015}, the design of hyperbolic metamaterials and perfect reflectors \cite{valentine2015} or the realization of zero-index metamaterials \cite{valentine2013}. It has also been demonstrated \cites{lippens2008,mirzaei2015} that the possibility of tailoring the artificial magnetic activity was a key to the design of invisibility cloaks based on dielectric materials. This has the great advantage of limiting the losses, as compared to metallic structures were strong losses are unavoidable, apart at the price of inserting active media \cite{Hess2013}.

The first mathematical study of this kind of dielectric metamaterials was made in the particular case where the structure is invariant in one direction
(see  \cites{prl,prcras} and \cite{random} for a generalization to the random case). In these papers
 the inclusions are infinite cylinders and the incident wave is polarized with a magnetic field parallel to the axis of the cylinders. 
As a consequence the original 3D-problem problem can be reduced in a two-dimension setting, this allowing a quite simple and rigorous asymptotic analysis.

\paragraph{Our contribution.}
In this paper we consider a dielectric structure placed in a \textit{bounded} domain of $\R^3$ 
and perform a full 3D asymptotic analysis. The infinitesimal parameter denoted $\eta$ represents the scale factor associated with the distance between inclusions. 
As it was announced in the note \cite{bou_bou_cras}, we will prove that the structure behaves, when $\eta\to0$, as a
\textit{local} material described by a frequency-dependent permeability tensor. 
Although this conclusion looks qualitatively in perfect agreement with what is obtained in the 2D case, the mathematical analysis reveals several novelties 
that we wish to emphasize here:

- In contrast with the 2D case,  the fast oscillations of the magnetic field are not anymore localized on the dielectric inclusions. 
 It follows that the induced magnetic activity cannot result simply from the superposition of independent Helmholtz micro-resonators as  
depicted in \cite{prl}. Interactions between the inclusions and the substrate have to be well understood. 

-  The averaging procedure we need for the asymptotic analysis has to be compatible with
the classical transmission conditions  across the boundary of the structure, namely the continuity of the tangential components of the electromagnetic field. 
This issue appears to be crucial regarding  the magnetic field: as it will be demonstrated, the bulk average (weak limit in $L^2_\mathrm{loc}$) will present tangential discontinuities. A key argument to overcome this difficulty will be the introduction of a new averaging process for periodic magnetic fields that we handle as one-forms on the complementary of the periodic inclusions.  
 
- The new spectral problem we have evidenced for describing the resonance modes of the structure is quite interesting by its interaction with the geometry.
 It involves divergence free periodic vector fields  on the three dimensional torus which are curl free on the complementary of the inclusions.
A direct numerical approximation problem of it turns out to be very costly computationally. A lot of attention has been devoted to finding
equivalent formulations making possible efficient numerical simulations for various type of shapes of dielectric inclusions. Some of these simulations are presented in Section~\ref{subsec:simul_mu}.

%
%
%

\paragraph{Notations.}\ 

\vspace{-0.1cm}
\begin{itemize}
 \item $\C^+:=\{z\in\C\ :\ \im(z)\ge0\}$,
 \item $\ov z$ denotes the complex conjugate of complex number  $z$, 
 \item $Y:=]-\frac12,\frac12[^3$,
 \item $\Sigma^*:=Y\setminus\ov\Sigma$,
 \item $B_R = \{x\in\R^3 \, :\, |x| <R\}$ \ ( $|x|$ denotes the Euclidean norm),
 \item $ A \Subset B$ for subsets of $\R^3$ means that $\ov A$ is a compact subset of the interior of $B$,
 \item $|B|$ denotes the Lebesgue measure of a Borel set $B\subset \R^3$,
 \item $1_B(x)$ denotes the characteristic function of $B$,
 \item $[\bfy]$ is the step function defined on $\R^3$ by $[\bfy]:=\bfk$ for all $\bfy\in Y+\bfk$, 
 \item $\langle u\rangle:=\int_Yu(y)\,dy$ denotes the mean value of a fonction $u\in L^1(Y)$,
 \item $\I$ is the identity matrix of order 3, 
 \item $\bs{M}:\bs{N}$ denotes the usual scalar product for $3\times 3$ matrices ,
 \item $\bfu\wedge\bfv$ denotes the cross product of vectors $\bfu,\bfv\in \R^3$,
 \item $\bfu\otimes\bfv$ denotes the tensor product of two vectors in $\R^3$,
\item $\disp\fint_B f(x)\,dx:=\frac{1}{|B|}\int_Bf(x)\,dx$ for all $f\in L^1_{\mathrm{loc}}(\R^3)$ and  Borel set $B\subset\R^3$. 
 
\item $C^\infty(\R^3)$ (resp.$C^\infty(K)$ if $K$ is a compact domain of $\R^3$),  denotes the space of functions which are $C^\infty$ on $\R^3$ (resp.on $K$) 
 \item $C^\infty_c(D)$ for an open subset $D\subset\R^3$, is the set of $C^\infty$-functions with compact support in $D$,
 \item $C^\infty_\sharp(Y)$ the subset of $Y$-periodic functions in $C^\infty(\R^3)$,
 \item $L^2_\sharp(Y):=\Big\{ f\in L^2_{loc}(\R^3)\ :\  f\ \text{$Y$-periodic}\Big\} $,
 \item $\langle \bff\cdot \bfg \rangle :=\ \int_Y \bff\cdot \ov{\bfg}\ dy\ $ the standard scalar product in $(L^2_\sharp(Y))^3$,
 \item $\Hunper(Y):=\left\{f\in L^2_\sharp(Y)\ :\ \partial_i f\in L^2_\sharp(Y)  \ ,\ i= 1,2,3\right\}$  \ ($\partial_i f$ is meant in the distributional sense).
 \item  $\Hunper(\Sigma^*)$ denotes the space of restrictions to $\Sigma^*$ of functions in $\Hunper(Y)$.
 \item  $\mathcal{D}'(A)$ denotes the distributions on the open subset $A\subset \R^3$.

 \end{itemize}

\paragraph{Geometrical assumptions.}
All along this paper the geometric domain of $\R^3$ in which the small dielectric inclusions are disposed will be denoted by $\Omega$.
This domain $\Omega$ is assumed to be bounded, simply-connected with Lipschitz boundary.
For every value of the small parameter $\eta>0$, we consider a diffracting obstacle occupying a subregion  $\Sigma_\eta \subset \dom$
which is obtained by periodization of a small inclusion of size $\eta$. 
 More precisely $\Sigma_\eta$ is given by
\begin{equation}\label{Ieta}
 \Sigma_\eta:=\bigcup_{k\in I_\eta}\eta(k+\Sigma)\quad,\quad I_\eta=\la i\in\Z^3\ |\  \eta(i+Y)\subset \dom\ra\ ,
\end{equation}
being $Y:=]-\frac{1}{2},\frac{1}{2}[^3$  the unit  cell and $\Sigma \subset\subset Y$ a reference inclusion. 

\med
The complexity of the diffracting obstacle is then encoded by the fast oscillating behavior of $\Sigma_\eta$ as $\eta$ becomes infinitesimal.
Let us notice that the \textit{filling ratio} of the inclusions remains positive when $\eta\to0$ since it converges to $|\Sigma|$ (the Lebesgue measure of $\Sigma$).

\begin{figure}
\begin{center}
\includegraphics[width=0.75\linewidth]{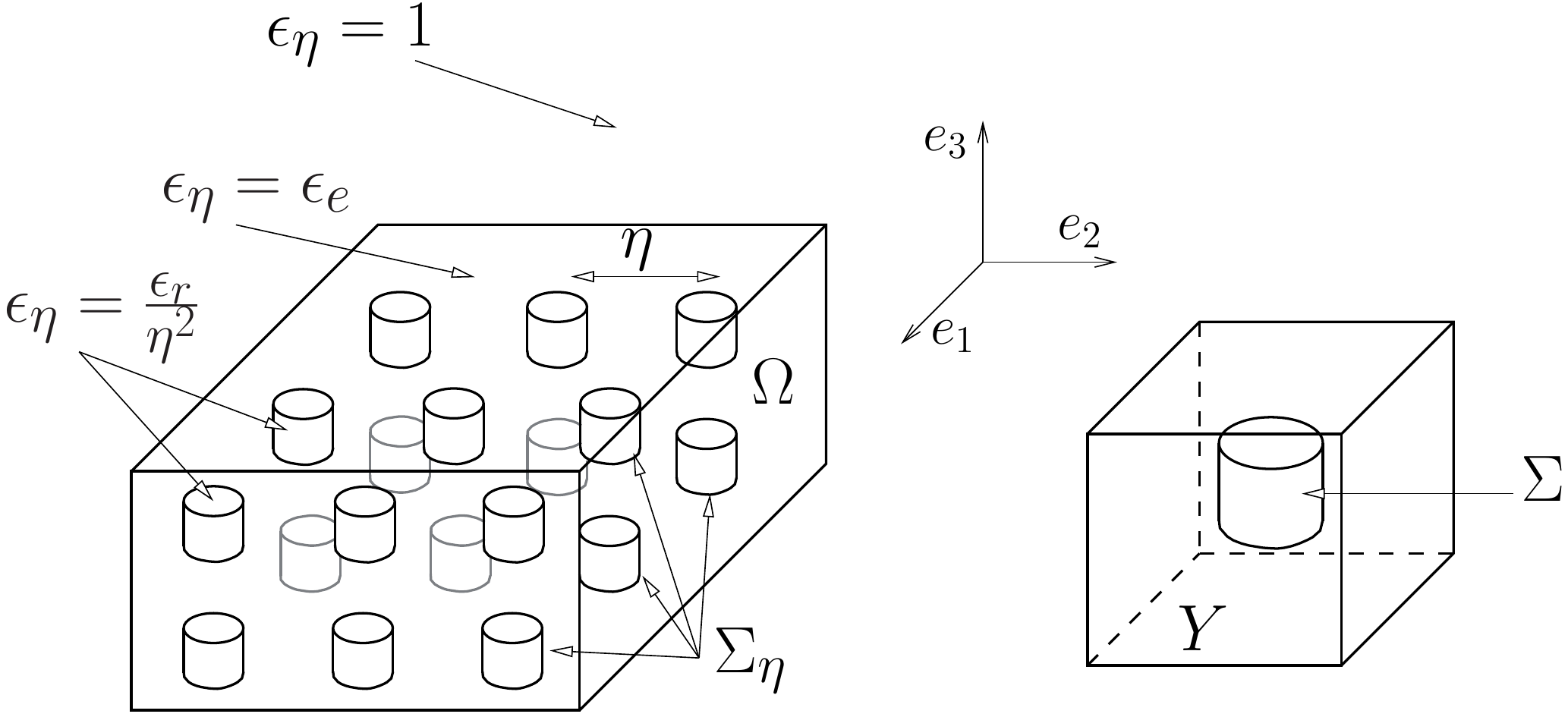}
\caption{Diffracting structure and  unit cell. \label{fig:omeg_mu}}
\end{center}
\end{figure}

\med
It turns out that the topology of the inclusion $\Sigma$  plays an important role in the asymptotic analysis as $\eta\to 0$. 
In this paper we will assume  that:
\begin{equation}\label{standingass}
\begin{array}{ll}
\text{ i)}  & \Sigma \ \text{is a connected compact subset of $Y$ with Lipschitz boundary} ,\\
\text{ii)}&  \Sigma^* = Y\setminus \Sigma\  \text{is simply connected.}
\end{array}
\end{equation}

\paragraph{Constitutive parameters and scaling.}

The local behavior of the medium is represented by its  relative permittivity and permeability tensors at every position $x\in \R^3$. In order to fit with the common use in the optical domain, we will assume the overall relative permeability to be constant equal to $1$. The dielectric properties of the structure under study are described by a function $\eps_\eta$ of the form   
\begin{equation}\label{def_epseta}
\eps_\eta(x):=\frac{\eps_r}{\eta^2}\,\un_{\Sigma_\eta}(x)+\eps_e\,\un_{\dom\setminus\Sigma_\eta}(x)+\un_{\R^3\setminus\dom}(x)\ ,
\end{equation} 
where the parameters $\eps_e\in\R^+$ and $\frac{\eps_r}{\eta^2} \in\C^+$ represent respectively the relative permittivity in the matrix and in the inclusions.
Here $\eps_r$ is a complex parameter such that:
\begin{equation}\label{hyp_im_eps_pos}
\Im( \eps_r)>0\ ,
\end{equation}
while the scaling factor $1/\eta^2$ is responsible of a high contrast becoming larger and larger as the period parameter $\eta$ of the structure decreases to zero.
The choice of this scaling is not new (see \cites{prl,prcras,random}). It ensures that the optical thickness of the inclusions remain constant and therefore the Mie resonances of each dielectric inclusion appear at frequencies which are independent of $\eta$ (see \cite{mie}).


\paragraph{Diffraction problem.}
The structure is illuminated by $(\bseinc,\bshinc)$ an incident monochromatic wave 
  travelling from infinity. We assume the harmonic time-dependence
to be $e^{-i\omega t}$ where $\omega>0$ is a fixed frequency. The total electromagnetic field
$(\bse_\eta,\bsh_\eta)$ satisfies Maxwell equations given in distributional sense in $\R^3$ by
\begin{equation} \label{prb:prin_mu}
\la
\begin{array}{l}
\rt \bse_\eta=i\omega\mu_0\bsh_\eta\ ,\\
\rt \bsh_\eta=-i\omega\eps_0 \, \eps_\eta \bse_\eta\ ,
\end{array}
\right.
\end{equation}
where $\eps_0>0$ and $\mu_0>0$ are  respectively the permittivity and permeability in the vacuum.
The influence of the incident wave is encoded by the fact that  the diffracted field $(\bse_\eta^d, \bsh_\eta^d):=(\bse_\eta-\bseinc$, $\bsh_\eta- \bshinc)$ 
satisfies the Silver-M\"uller's condition at infinity
\begin{equation} \label{SM}
(\bse_\eta^d, \bsh_\eta^d)   =   O\left(\frac{1}{|\bsx|}\right), \qquad
\omega \eps_0 \left( \frac{\bsx}{|\bsx|}\wedge \bse_\eta^d \right) - k_0 \bsh_\eta^d
 =  o\left(\frac{1}{|\bsx|}\right).
\end{equation}

\begin{remark}\label{rem:adim_mu}\nshape{
In this paper we will proceed in a dimensionless framework assuming implicitly that the physical period of the composite is in fact $\eta\, d$
being  $d$ the unit of length. Thereby the adimensional wavelength of the electromagnetic waves becomes $\lambda/d$ (if $\lambda$ is the real wavelength).
In order to simplify notations this parameter $d$ will not appear in the following, except in the presentation of numerical simulations in Section \ref{subsec:simul_mu}.
}\end{remark}  

%

\section{Presentation of the results}\label{sec:reslut_prin_mu}


The asymptotic analysis as $\eta\to 0$ of $(\bse_\eta,\bsh_\eta)$ solving \eqref{prb:prin_mu} 
leads to a homogenized diffraction problem of the kind
\begin{equation}\label{prob_lim_mu_intro}
\la
\begin{array}{l}
\rt \bs{E}=i\omega\mu_0\,\bs{\mu}(x,\omega)\,\bs{H}\\
\rt \bs{H}=-i\omega\eps_0\,\bs{\eps}(x)\,\bs{E}\\
(\bse-\bseinc,\bsh-\bshinc)\mbox{ satisfies  \eqref{SM}}
\end{array}
\right.
\end{equation}
where
\begin{equation}\label{epsmux}
\bs{\eps}(x):=\bs{I}_3\,1_{\R^3\setminus\dom}(x)+\epsef\,1_{\dom}(x)\quad,\quad \bs{\mu}(x,\omega):=\bs{I}_3\,1_{\R^3\setminus\dom}(x)+\muef(\omega)\,1_{\dom}(x)\ ,
\end{equation} 
Here the effective tensors $\epsef, \muef$ describe for each frequency $\omega$ a homogeneous medium occupying the domain $\dom$.
They are described in a precise way in the next subsection. The convergence of $(\bse_\eta,\bsh_\eta)$ to the solution of \ref{prob_lim_mu_intro}
will be specified in our main Theorem (Theorem \ref{t.main}).
 
Let us point out that the two first equations in \eqref{prob_lim_mu_intro} are to be understood in the distributional sense in $\R^3$.
In particular, under mild regularity assumptions, they imply the following transmission conditions on $\partial\dom$:
\begin{equation}\label{transmission}
[ \bfn\wedge \bse] = [ \bfn\wedge \bsh] = 0 \ ,\qquad [ \bfn\cdot\bmu\bsh]=[ \bfn\cdot\beps\bse]=0 \ ,
\end{equation}
with $\bfn$ denoting the outward unit vector and $[\cdot]$ the jump across $\partial\dom$.

\med
The proof of the
existence and uniqueness of the solution to \eqref{prb:prin_mu}, \eqref{SM}  
under the dissipativity condition \eqref{hyp_im_eps_pos} is classical and can be found
e.g. in \cite{cessenat}.  With regard to the uniqueness for a limit problem of the form
given in \eqref{prob_lim_mu_intro}, we have
\begin{lemma}\label{wellposed}
Assume that $\epsef$ is real symmetric positive and that  $\muef$ is a symmetric tensor whose imaginary part $ \Im(\muef)$ is positive definite. Then
the solution to \eqref{prob_lim_mu_intro} is unique. 
\end{lemma}
\proof
By linearity, it is enough to check that if  $(\bse,\bsh)$ solves \eqref{prob_lim_mu_intro} for a vanishing $(\bseinc,\bshinc)$, then  $(\bse,\bsh)=(0,0)$.
Let $R>0$ so large that $\dom\subset \subset B_R$ and denote  $\mP(R):= \int_{\partial B_R} (\bse\wedge\ov\bsh)\cdot\bfn $ the flux of the Poynting vector. 
As usual, the real part of $\mP(R)$ does not depend on $R$. Indeed exploiting the identity  $\dv \bse\wedge\ov\bsh = \rt\bse\cdot\ov\bsh-\rt\ov\bsh\cdot\bse$,
 we may integrate by parts and , by taking into account \eqref{prob_lim_mu_intro}, we get 
\begin{equation}\label{balance}
\mP(R)\ =\ \int_{B_R}\left(\rt\bse\cdot\ov\bsh-\rt\ov\bsh\cdot\bse\right)\ =\ i\omega\Big(\mu_0 \int_{B_R}\bmu\bsh\cdot\ov\bsh-\eps_0\int_{B_R}\ov \beps\bse\cdot\ov\bse\Big)\ .
\end{equation}
In particular, since $\beps$ is real and $\bmu$ agrees with the identity tensor outside $\dom$, by identifying the real and imaginary parts, we deduce that
\begin{equation}\label{pointing_unicite}
\Re\lp \mP(R) \rp\ =\ -\omega\mu_0\int_\dom\Im\Big(\muef\bsh\cdot\ov\bsh\Big)\ .
\end{equation}
We may now  pass to the limit $R\to +\infty$ in the left-hand member of \eqref{pointing_unicite}. Now by exploiting the fact that $(\bse,\bsh)$
satisfies \eqref{SM} with $\Einc=\Hinc=0$, we  find that $\lim_{R\to\infty} \mP(R) = 0$.  Thus  the left-hand member of \eqref{pointing_unicite} vanishes
as well as the integral of $\Im\Big(\muef\bsh\cdot\ov\bsh\Big)$ over $\dom$. By the positivity assumption on $\Im(\muef)$ and since $\muef$ is symmetric, we have
for a suitable constant $c>0$:
$$  \Im\Big(\muef\bsh\cdot\ov\bsh\Big) \ =\ \Im(\muef) \bsh\cdot\ov\bsh \ \ge\  c\, |\bsh|^2 \ ,$$
 Thus $\bsh$ vanishes on $\dom $. By the second equation of \eqref{prob_lim_mu_intro}
 and the fact that $\beps$ is a real positive tensor, it is also the case of $\bse$. It is then classical to deduce that $(\bse,\bsh)$ vanishes in the exterior domain $\R^3\setminus \dom$ as well.
\qed

\subsection{Effective laws.}

The limit diffraction problem we wrote in the form \eqref{prb:prin_mu} is completely determined by relations \eqref{epsmux} once we know the 
 effective permittivity tensor $\epsef$ and the
effective permeability tensor $\muef$. These tensors  are described in a precise way in the next two paragraphs. 
It turns out that $\epsef$ is real positive and \textit{does not depend} on the frequency and on the dielectric parameter $\eps_r$. 
In contrast, the tensor $\muef= \muef(\omega)$ depends on the frequency and exhibits resonances.

\paragraph{Effective permittivity law.}  It depends only on the geometry of $\Sigma$ and on the permittivity in the matrix surrounding inclusions. The computation of $\epsef$ 
looks similar as the one used in the two-dimensional case (see \cites{prcras,random}) where the classical ingredients of homogenization theory
for Neumann problems with holes can be recognized. The entries of the tensor $\epsef$ are given for $(k,l)\in\{1,2,3\}^2$ by
\begin{equation}\label{def_epsef_intro}
\epsef_{kl}:=\eps_e\int_Y(e_k+\nabla\chi_k).(e_l+\nabla\chi_l)\ ,
\end{equation}
being $\chi_k\in W^{1,2}_\sharp(Y;\R)$ the unique solutions of
\begin{equation}\label{localE}
\Delta_y\chi_k=0\quad\mbox{in}\quad \Sigma^*\qquad\mbox{and}\qquad \chi_k=-y_k\quad\mbox{in}\quad \Sigma\ .
\end{equation}
As $\eps_e$ is a positive real, it can be readily checked (see \eqref{positive-eps}) that the tensor  $\epsef$ is  real symmetric positive.

\paragraph{ Effective permeability law.} The dependence of permeability tensor $\muef$ with respect to $\omega$ is ruled by the internal resonances of the composite structure which are responsible for the magnetic activity. The description of the underlying spectral problem 
is quite involved due to the fact that strong oscillations of the microscopic magnetic field $\bsh_0(x,\cdot)$ are allowed not only in $\Sigma$ (in a similar way as in the 2D case \cites{random, prcras, prl}),
but also in the surrounding matrix. A nice way to circumvent this difficulty consists in looking at the curl of the magnetic field which accounts for the magnetic activity.
It turns out that this curl vanishes outside $\Sigma$ (see \eqref{PH02}). We therefore introduce the space $Z_0\subset L^2(Y;\R^3)$ defined by
\begin{equation}\label{def_Z0}
 Z_0:=\Big\{\bs{f}\in L^2(Y;\R^3)\ :\ \dv \bs{f}=0,\quad \bs{f}=0 \quad \mbox{in $\Sigma^*$}\Big\}\ .
\end{equation}
Let us notice that the elements of $Z_0$ can be identified with divergence-free vector fields in $L^2(\Sigma;\R^3)$
with vanishing  normal trace on $\partial\Sigma$.
Next, we associate to every element $\bff\in Z_0$  the unique solution $\bpsi_f$ in $\Hunper(Y;\R^3)$ of
\begin{equation}\label{def:phi_f}
-\Delta \bpsi_f=\bff \quad \mbox{in }Y\ ,\quad
\int_Y\bpsi_f=0\ .
\end{equation}
Then, as will be discovered later, the resonance frequencies for the microscopic magnetic field are directly related to
 the following eigenvalue  problem. Find $(\bff,\alpha)\in Z_0\times\R$ such that for all $\bfg\in Z_0$:
 \begin{equation}\label{pos_A}
\int_Y\nabla\bpsi_f:\nabla\bpsi_g \ +\ \frac{1}{4}\ \lp \int_\Sigma \bfz\wedge\bff \,d\bs{z}\rp\cdot\lp \int_\Sigma \bfz\wedge\bfg \,d\bs{z}\rp\ = \ \alpha \, \int_\Sigma \bff\cdot\bfg
\end{equation}
The linear operator associated with the bilinear form in the left-hand side  turns 
out to be positive, compact and self-adjoint on the Hilbert space $Z_0$ (embedded with the $L^2(Y;\C^3)$ scalar product).
Therefore, it exists a sequence of eigenvalues $\alpha_0 \ge \alpha_1 \ge \dots \ge \alpha_n \ge \dots >0$ such that $\alpha_n \to 0$ and an associated orthonormal basis of eigenvectors 
$\{\bff_n,\ n\in\N\}$ in $Z_0$. 

The effective permeability law we are going to establish for the limit diffraction problem is described by a symmetric tensor $\muef$. This tensor
can be then written as the following series:
\begin{equation}\label{def_muef_intro2}
\muef(k_0)\ :=\ \I+ \frac1{4}\, \sum_{n\in \N}\frac{\eps_r k_0^2}{1-\eps_r\alpha_n k_0^2}\lp\int_\Sigma \bs{y}\wedge\bff_n\,dy\rp\otimes\lp\int_\Sigma \bs{y}\wedge\bff_n\,dy\rp \ .
\end{equation}
In fact, it is convenient to present an alternative representation of $\muef$ involving periodic vector fields $\bfu_n$ on the unit cell which will be useful to describe the
fast oscillating magnetic field $\bsh_\eta$. Let us define
\begin{equation}\label{wn}
\lambda_n= \frac1{\alpha_n}\quad,\quad \bfu_n = \frac1{\sqrt{\alpha_n}}\Big( \rt(\bpsi_{f_n}) + \frac1{2} \int_\Sigma \bfy\wedge \bff_n \Big)\ .
\end{equation}
Then the following relation holds
\begin{equation}\label{def_muef_intro}
\muef(k_0)\ =\ \disp\I+\sum_{n\in \N}\frac{\eps_rk_0^2}{\lambda_n-\eps_rk_0^2}\lp\int_Y\bfu_n\rp\otimes\lp\int_Y\bfu_n\rp\ .
\end{equation}
As will be seen later in Section \ref{sec:new_form_spec_num}, the pair $(\bfu_n,\lambda_n)$ can be characterized directly as solutions of  the 
following spectral problem. Find $(\bfw,\lambda)\in X_0^{\dv}\times\R$ such that for all $\bfv\in X_0^{\dv}$:
\begin{equation}\label{spectralw}
\int_{Y}\rt\bfu_n\cdot\rt\bfv=\lambda_n\int_Y\bfu_n\cdot\bfv\ ,
\end{equation}
where $X_0^{\dv}$ is a suitable subspace of $\Hunper(Y;\C^3)$ consisting of functions which are curl-free in $\Sigma^*$ and divergence-free in $Y$ (see \eqref{def_X}).

\subsection{Main convergence result.}

\med In view of Lemma \ref{wellposed}, we consider the unique solution $(\bse,\bsh)$ of \eqref{prob_lim_mu_intro}  and  denote by $E_k,H_k$  the $k$-th component of 
$\bse,\bsh$ respectively ($k\in\{1,2,3\}$). Next, we introduce two important vector fields $\bse_0(x,y)$, $\bsh_0(x,y)$ in $L^2(B_R\x Y)$ (associated with the
two-scale analysis performed in Section \ref{sec:esti_prelim_mu}) where $x$ represents the macroscopic variable and where a $Y$-periodic dependence with respect to the fast variable $y$ is set in order
to account for the oscillating behavior of the sequence $(\bse_\eta,\bsh_\eta)$.

\med The ``two-scale electric field'' is defined by
\begin{equation}\label{E0}
\bse_0(x,y)= \begin{cases} \disp \ \sum_{k=1}^3E_k(x)\big( \bfe_k+\nabla\chi_k(y) \big)\ & \text{if $x\in \dom$} \\
\ \bse(x) & \text{if $x\in B_R\setminus \dom$}
\end{cases} \end{equation}
where the functions $\chi_k$ are the solutions of \eqref{localE}.
Similarly, with the help of the periodic vector fields $\bfu_n$  and positive numbers $\lambda_n$  defined in \eqref{wn}, 
we define  the ``two-scale magnetic field''  by 
\begin{equation}\label{H0}
\bsh_0(x,y)=
\begin{cases} \disp \ \sum_{k=1}^3H_k(x)\, \bsh^k(y)\  & \text{if $x\in \dom$} \\
\ \bsh(x) & \text{if $x\in B_R\setminus \dom$}
\end{cases}
\end{equation}
\begin{equation}\label{Hk} \bsh^k(y) := 
\bfe_k+\sum_{n\in\N} \langle\bfe_k,\bfu_n\rangle\frac{\eps_r k_0^2}{\lambda_n-\eps_r k_0^2}\bfu_n(y)\end{equation}

\noindent
We are now in a position to state the main result of the paper:
\begin{theorem} \label{t.main}\ \\
Let us assume \eqref{standingass}, \eqref{hyp_im_eps_pos} and let $\bs{\eps}, \bs{\mu} $ be defined by  \eqref{epsmux}, \eqref{def_epsef_intro}, 
\eqref{def_muef_intro}. Let $(\bse,\bsh)$ be the unique solution of \eqref{prob_lim_mu_intro}. Then  the solution $(\bse_\eta,\bsh_\eta)$ of the diffraction problem \eqref{prb:prin_mu} 
satisfies
\begin{equation}\label{strong_conv}
\int_{B_R}\Big| \bsh_\eta(x)-\bsh_0\xxeta \Big|^2\dx\to0\ ,\qquad
\int_{B_R}\Big|\bse_\eta(x)-\bse_0\xxeta \Big|^2\, dx\to0\ .
\end{equation}
where $R$  is arbitrary large and $\bse_0, \bsh_0$ are given by \eqref{E0} and \eqref{H0} respectively.
Futhermore it holds  $(\bse_\eta,\bsh_\eta)\to (\bse,\bsh)$ in $C^\infty(K)$ for every compact subset  $K\Subset \R^3\setminus\dom$.  
\end{theorem}

The proof of Theorem \ref{t.main} is quite long and involved. It is postponed to Section \ref{sec:proof_prin} where the arguments are presented along two steps.
The most delicate issue is the $L^2$ upper-bound estimate (see \eqref{born_prio1}) for the electromagnetic field. It is proved \textit{a posteriori} in the last step by using a contradiction argument (in the same line as in \cites{boufel_cap,bou_schw}). 
Before this proof, in Section \ref{sec:electromagnetic},  we assume \textit{a priori} this $L^2$- upper-bound in order to prepare the two-scale analysis of the system.

\begin{remark}\label{bulkaverage} \nshape{ Let us emphasize that the convergence result in Theorem \ref{t.main} is unusual in the classical framework of homogenization theory:
the effective magnetic field $\bsh$ that we use in order to describe the limiting diffraction problem \eqref{prb:prin_mu}  does not agree inside the obstacle with the weak limit of $\bsh_\eta$ in $L^2_{\mathrm{loc}}(\R^3)$. 
Indeed in view of  \eqref{def_muef_intro}, \eqref{H0} and \eqref{Hk}, it is easy to check that, for $x\in \dom$, this weak limit $\tilde \bsh$ satisfies:
\begin{equation} \label{bulkmean}
{\tilde \bsh}(x) \ :=\ \int_Y\bsh_0(x,y)\,dy\ =\ \muef\,\bsh(x)\ , 
\end{equation}
whereas $\tilde \bsh= \bsh$ in $\R^3\setminus \Omega$. As shown later, this  tensor $\muef $ differs from the identity matrix (for most of the frequencies). 

\med There is a major reason  in not using $\tilde \bsh$ for describing the limit magnetic field: 
 the tangential trace of $\tilde\bsh$ on $\partial\dom$ turns out to differ
from that of the field $\bsh$ outside.  This is {\it a priori} not physically reasonable and suggests
 that taking the asymptotic bulk average of $\bsh_\eta$ would not be a good choice.
  Moreover, in view of Theorem \ref{t.main} and of equations \eqref{prob_lim_mu_intro} solved by $(\bse,\bsh)$,
   we find \textit{a posteriori} that the limit system written in term of $\tilde \bsh$ leads to
$$
\rt \bs{E}=i\omega\mu_0\,\tilde\bsh\ ,\qquad
\rt (\bmu^{-1}(x) \, \tilde\bsh) = -i\omega\eps_0\,\bs{\eps}(x)\,\bs{E}\ ,
$$
where we lose the $\rt$ structure of the second Maxwell and the fact that the magnetic activity is encoded through the tensor $\muef$.

\med
 In order to obtain a proper notion of effective magnetic field, we will use a different averaging recipe in which
the periodic field $\bsh_0(x,\cdot)$ is seen as a  closed periodic differential 1-form on $\Y\setminus\Sigma$.  This allows to define at every point $x\in\dom$  a circulation vector $\bsh(x) =\oint \bsh_0(x,\cdot)$ (see Lemma \ref{lem:def_circulation}). Adopting this alternative definition of $\bsh$ inside $\dom$, we will succeed in recovering the classical transmission conditions across $\partial\dom$  namely \eqref{transmission}.
  We notice that in contrast, the limit \textit{magnetic induction} vector field $\bsB := \mu_0 \muef \bsh$ agrees with the weak limit of $\bsB_\eta = \mu_0 \bsh_\eta$. Roughly speaking $ \bsB_\eta$ can be seen as a closed $2$-form and as $\eta\to 0$ the  local flux it generates is represented by the bulk average $\bsB$.

}\end{remark}
 
%
%
%

\medskip
\subsection{Frequency-dependent permeability and band gaps.}

According to Theorem \ref{t.main}, for infinitesimal $\eta$, the electromagnetic field $(\bse_\eta,\bsh_\eta)$ outside the obstacle is close to the solution of a limit diffraction problem in which domain $\dom$ is occupied by a homogeneous medium whose
permeability and permittivity tensors are given  in \eqref{def_epsef_intro} and  \eqref{def_muef_intro} respectively. 
In this asymptotic model, the most interesting issue with respect to applications stems from the properties of  tensor $\muef$,
in particular its explicit dependence with respect to  the angular frequency  $\omega=k_0(\eps_0 \mu_0)^{-1/2}\ $
as well as its ability to exhibit eigenvalues with a \textit{negative real part}. To see that, it is convenient 
to introduce $V_\lambda$ the eigenspace associated with a  eigenvalue
$\lambda$ of spectral problem  \eqref{spectralw}. Then denoting by  $P_{V_\lambda}$  the 
orthogonal projector on  $V_\lambda$ (with respect to the scalar product of $L^2(Y)^3$),
we may rewrite \eqref{def_muef_intro} as follows
\begin{equation}\label{mu_vp}
  \muef(\omega):=1+\sum_{\lambda\in\sigma_0}\frac{\eps_rk_0^2}{\lambda-\eps_rk_0^2}M_\lambda\quad , \quad (M_\lambda)_{kl}:=(P_{V_\lambda}(\bfe_k),\bfe_l)\ ,
\end{equation}
where 
\begin{equation}\label{sigma0}
\sigma_0 :=\left\{ \lambda \ :\quad  \lambda \text{ eigenvalue of (\ref{spectralw})} \quad ,\quad  M_\lambda\not=0 \right\} 
\end{equation}
For $\lambda\in \sigma_0$, the real symmetric matrix satisfies $0\le M_\lambda\le I_3$. It is of rank one if $\lambda$ is simple. In order for it to be of full rank
 we need at least that the multiplicity of $\lambda$ be not smaller than $3$.

Let us assume for simplicity that  
 $\eps_r$ is a positive real (lossless dielectric inclusions). Then tensor $\muef(\omega)$ is  real symmetric, continuous with respect to $\omega$
except at the frequencies \textit{ (Mie resonances)} given by:  
$$\omega_\lambda:=\sqrt{\lambda}
(\eps_0\mu_0\eps_r)^{-1/2}\quad,\quad \lambda\in\sigma_0\ .$$ 
In the vicinity of this values, $\muef$ blows up and  we are led to different  consequences according to the rank of $M_\lambda$.
Denote by  $\mu^\pm$ the largest (resp. the smallest) of the  eigenvalues of $\muef$. Clearly it follows from \eqref{mu_vp} that
$
\lim_{\omega\to\omega_\lambda^+}\mu^-(\omega)=-\infty 
$. The same holds true 
for $\mu^+$ if and only if  $M_\lambda$ has full rank.
In this case, by continuity, we obtain an interval of frequencies in which all eigenvalues of $\muef$ are negative.
For such frequencies, since $\epsef$ is  positive definite, the electromagnetic field cannot propagate in any direction inside the obstacle. We may therefore conclude to the existence of \textit{ a photonic band gap}.
On the opposite side if  $M_\lambda$ is not of full rank, then vectors in its kernel determine propagative directions for the electromagnetic field.
Such a partial band gap situation was already observed in the context of elastic waves \cite{miara}.

\begin{figure}[ht]
\includegraphics[trim=1.7cm 0cm 2cm 0.8cm,clip,width=0.49\linewidth]{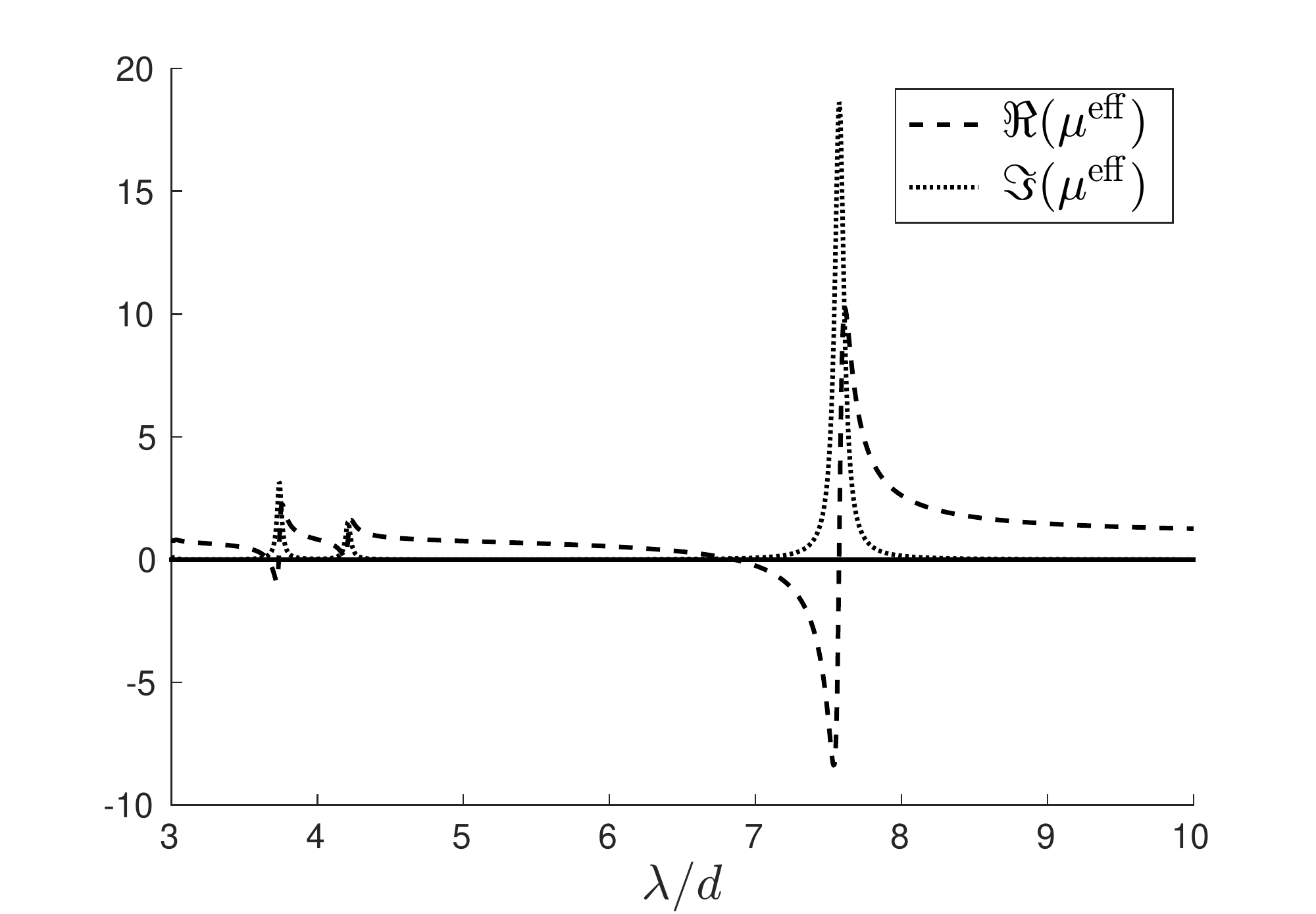}
\hfill
 \includegraphics[trim=1.7cm 0cm 2cm 0.8cm, clip,width=0.49\linewidth]{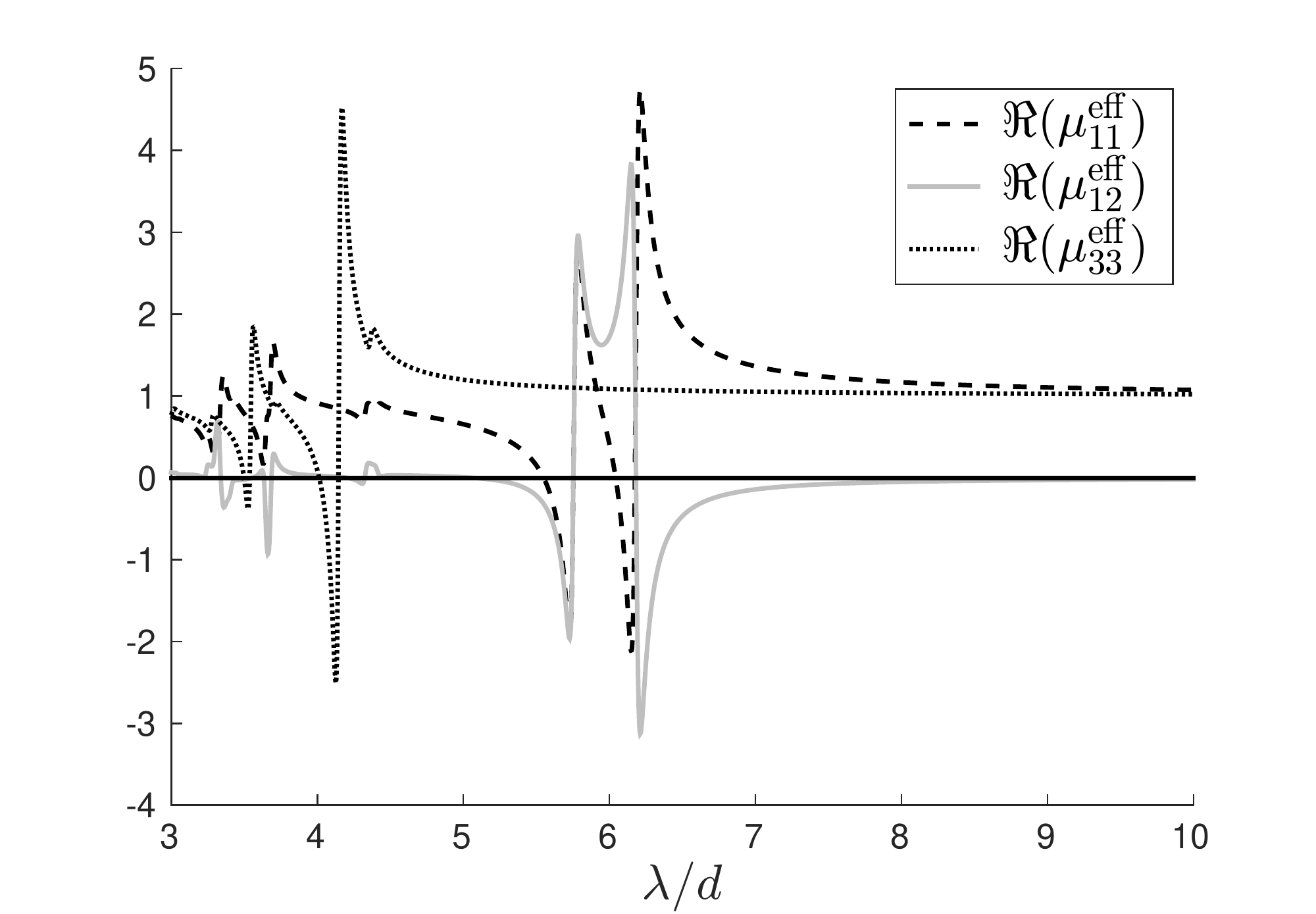}
 
\caption{\small Tensor $\muef$ as a function of $\lambda/d$ for $\eps_r=100+i$.
On the left  $\Sigma=\Sigma_1$ is a cube of size $0.6$ and tensor $\muef$ is scalar; on the right $\Sigma=\Sigma_2$ is a L-shape inclusion (and
$\muef_{22}=\muef_{11}, \muef_{13}=\muef_{23}=0$).  \label{fig:mu_1}}
\end{figure}

\section{Numerical simulations.}\label{subsec:simul_mu}

In this Section we present a numerical approach in order to evaluate the tensor $\muef$ as a function of the frequency.
To that aim it is convenient to use the representation  \eqref{def_muef_intro2} where we need to solve the three dimensional cell spectral problem \eqref{pos_A} 
in the space $Z_0$ of divergence-free fields vanishing outside $\Sigma$. The advantage to work with this representation rather than with \eqref{def_muef_intro}
is that the space $Z_0$ requires a discretization on subset $\Sigma$ only. 
The  approximation of \eqref{pos_A} is performed by means of a Galerkin method making use of the piecewise affine edge-elements of Nedelec (see \cite{nedelec}). 
The main drawback of the method  with respect to the computation cost is that we need to handle the non-local 3D-cell problem \eqref{def:phi_f}. For this problem we
 use an integral equation method with the help of the  Green kernel of the  inverse Laplace operator on the $3$-dimensional torus (we used the explicit form 
given in  \cite{Marshall}).

In order to give a nice description of $\muef(\omega)$ in the vicinity of the resonant frequencies, we have computed some of the local displacement currents $\bsj^k,
1\le k\le 3$ appearing in the periodic cell problem (see \eqref{decomp_Hk}). Recall that $\bsj^k$ represents the vorticity of the shape magnetic field $\bsh^k$ and that it belongs to space $Z_0$. By solving the spectral problem  \eqref{pos_A}, we obtain a sequence of pairs $(\bff_n,\alpha_n)$ and then we recover $\bsj^k= \rt \bsh^k$
from the following expansion:
 \begin{equation}\label{Jk}
 \bsj^k\ =\ \frac12\sum_{n\in\N} \frac{\eps_r\,k_0^2}{1-\alpha_n\,\eps_r\,k_0^2}\langle \bfe_k,\bfy\wedge\bff_n \rangle\,\bff_n\ .
\end{equation}
The equality above is a straightforward consequence of \eqref{Hk} taking into account that,  by \eqref{wn}, one has 
$$ \rt \bsu_n = \frac{\bff_n}{\sqrt{\alpha_n}}\quad ,\quad \langle \bfe_k, \bsu_n\rangle = \frac1{2 \sqrt{\alpha_n}} \langle \bfe_k,\bfy\wedge\bff_n \rangle\ .$$ 
Then it follows from   \eqref{def_muef_intro2} and \eqref{Jk} that  the entries of tensor $\muef(\omega)$ can be deduced from the relations
\begin{equation}\label{muefJ}
\muef_{kl}=\frac12\langle \bfy\wedge\bsj^k,\bfe_l\rangle\qquad  k,l \in \{1,2,3\}\ .
\end{equation}
For a given frequency, the excited resonances are weighted by the strength factor $\int_\Sigma\bfy\wedge\bff_n$ so that, in practice, only a few of the eigenvectors $\bff_n$ will be contributing
in  the expansion \eqref{Jk}. We are going to represent them showing thus the significative loops of displacement current $\bsj^k$ which take place in each periodic cell of the composite structure.

\med
In view of  \eqref{Jk} and \eqref{muefJ}, numerical simulations for solving \eqref{pos_A} 
have been performed for two kind of geometries of the dielectric inclusion $\Sigma$ that we will discuss separately: 
\begin{enumerate}
\item[-] 
 the cubic one where \quad  $\Sigma =\Sigma_1:=(-0.3,0.3)^3$.
\item[-] the $L$-shaped one where \  $\Sigma=\Sigma_2:= (-0.3,0.3)^3\setminus \big([-0.3,0.1]^2\times\R\big)$.
\end{enumerate}
In both cases the permittivity parameter $\eps_r$ characterizing the dielectric inclusion is taken to be  $100+i$, meaning that
the structure is slightly \textit{dissipative}.

\begin{figure}[!h]
 \includegraphics[trim=17cm 3.5cm 15cm 2cm,clip,width=0.49\linewidth]{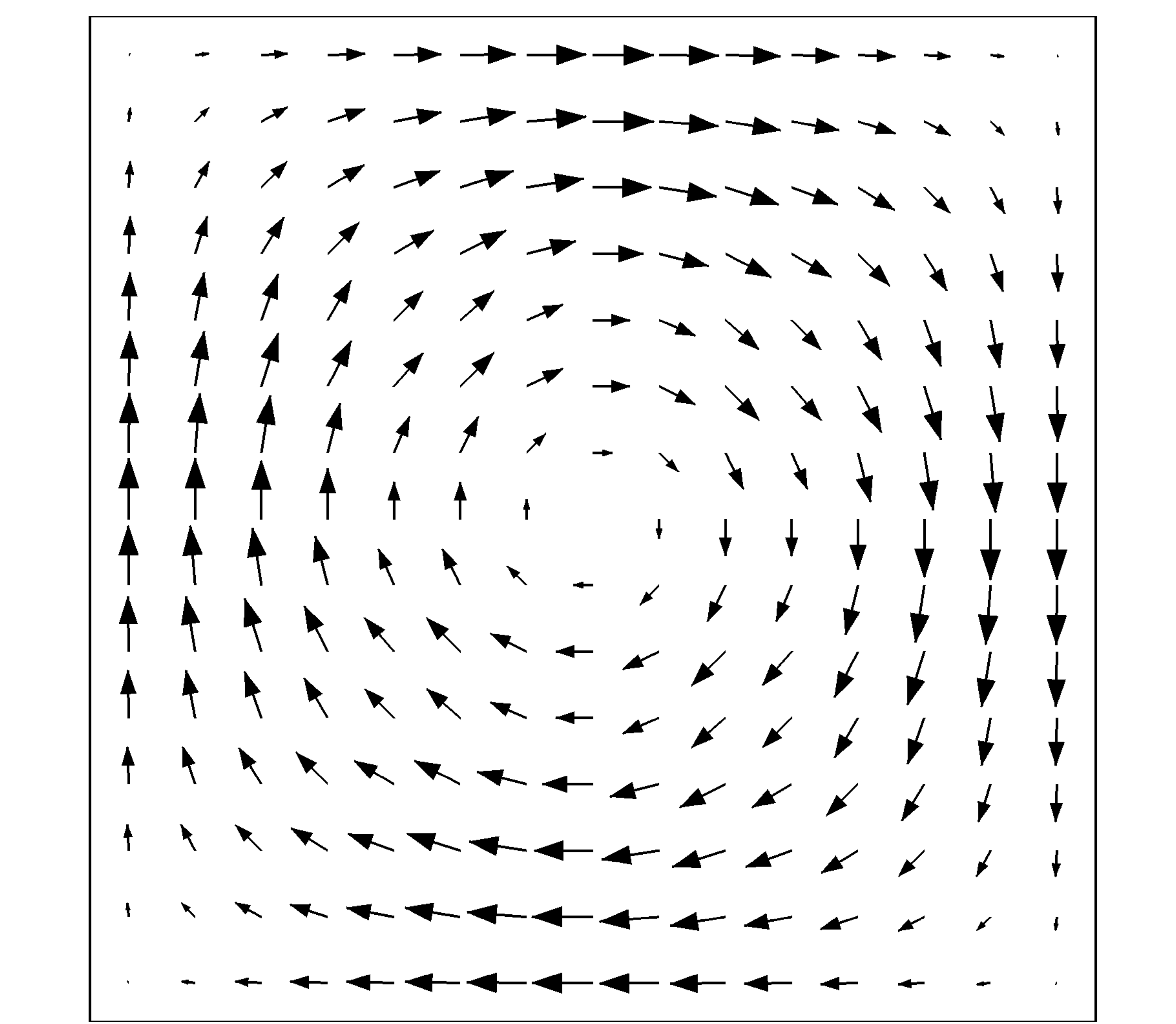}
 \hfill
\includegraphics[trim=17cm 3.5cm 15cm 2cm,clip,width=0.49\linewidth]{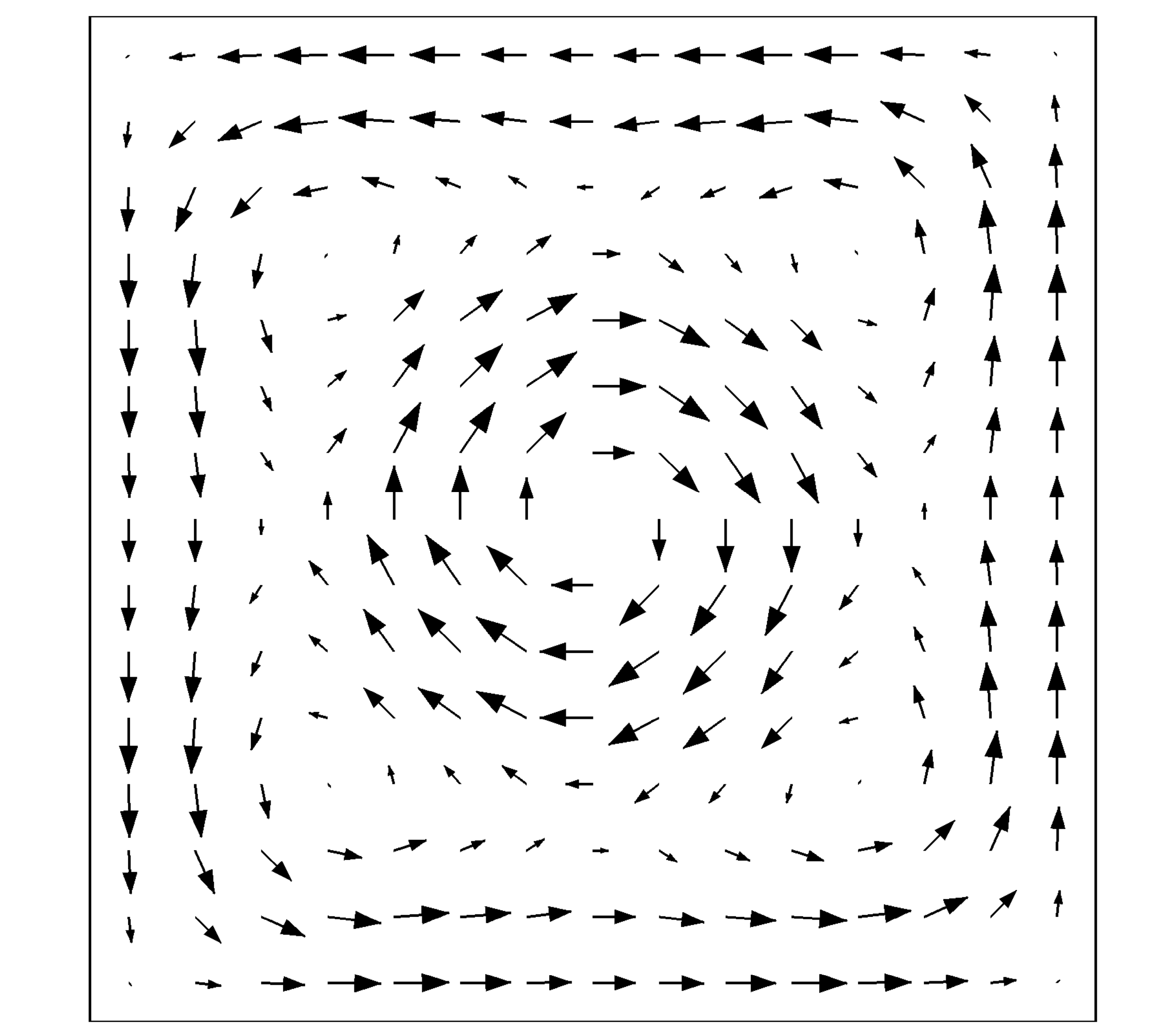}
 \caption{\small \ Horizontal section of $\bsj^3$ for the cubic geometry. Left side: $\lambda/d=7.5$ (fundamental resonance). Right side:
$\lambda/d=3.76.$ (resonance associated with the $24$th eigenvalue).  \label{fig.sectionJ3}}
\end{figure}

\paragraph{Cubic geometry.} As  the faces of the cubic inclusion $\Sigma_1$ and of the periodic cell $Y$  share the same orientation,
many symmetry properties can be exploited, namely the invariance of the structure under all rotations of angle $\pm\pi/2$ with axis $\bfe_i, i\in \{1,2,3\}$.
In particular it is easy to deduce that $\muef (= \muef I_3)$ is a scalar tensor. Moreover  an eigenvalue $\alpha$ such that
 $\int_\Sigma\bfy\wedge\bff \neq 0$ for any associated eigenvector $\bff$
has a multiplicity $\ge 3$. If it is not the case, that means that $\alpha $ is not contributing in expansion \eqref{def_muef_intro2}
 (equivalently $\lambda=\alpha^{-1}$ does not belong to $\sigma_0$ defined in \eqref{sigma0}).
In figure \ref{fig:mu_1}, we represent on the left side the real and the imaginary parts of $\muef$ as  a function of the normalized wavelength $\frac{\lambda}{d}$ (see Remark \ref{rem:adim_mu}).
Band gaps correspond to the frequency intervals in which the scalar permeability $\muef$ satisfies $\Re(\muef)<0$.
We represent only the part of the graph where the influence of resonant frequencies is significant namely the part corresponding to $\frac{\lambda}{d} \in [3,10]$ 
in which three oscillations of $\muef$ can be observed. Two of them have a sufficiently large amplitude in order to reach a  negative $\Re(\muef)$:
the larger one corresponds to the fundamental eigenvalue $\alpha_1$ while the second one is associated with eigenvalue $\alpha_{24}$
(here the $\alpha_n$ are repeated accounting their multiplicity). The third one associated with $\alpha_{17}$ corresponds to a resonance whose amplitude is too small
to force an additional change of sign for $\Re(\muef)$. In fact the numerical computations reveal that a very few of the $\alpha_n$ contribute to the series: among the $49$ first ones
only $\{\alpha_1,\alpha_2,\alpha_3,\alpha_{17},\alpha_{18},\alpha_{19},\alpha_{24},\alpha_{25},\alpha_{26}\}$ are such that $\alpha^{-1} \in \sigma_0$.

\begin{figure}[ht]
 \includegraphics[trim=6cm 5cm 6cm 4cm, clip,width=0.49\linewidth]{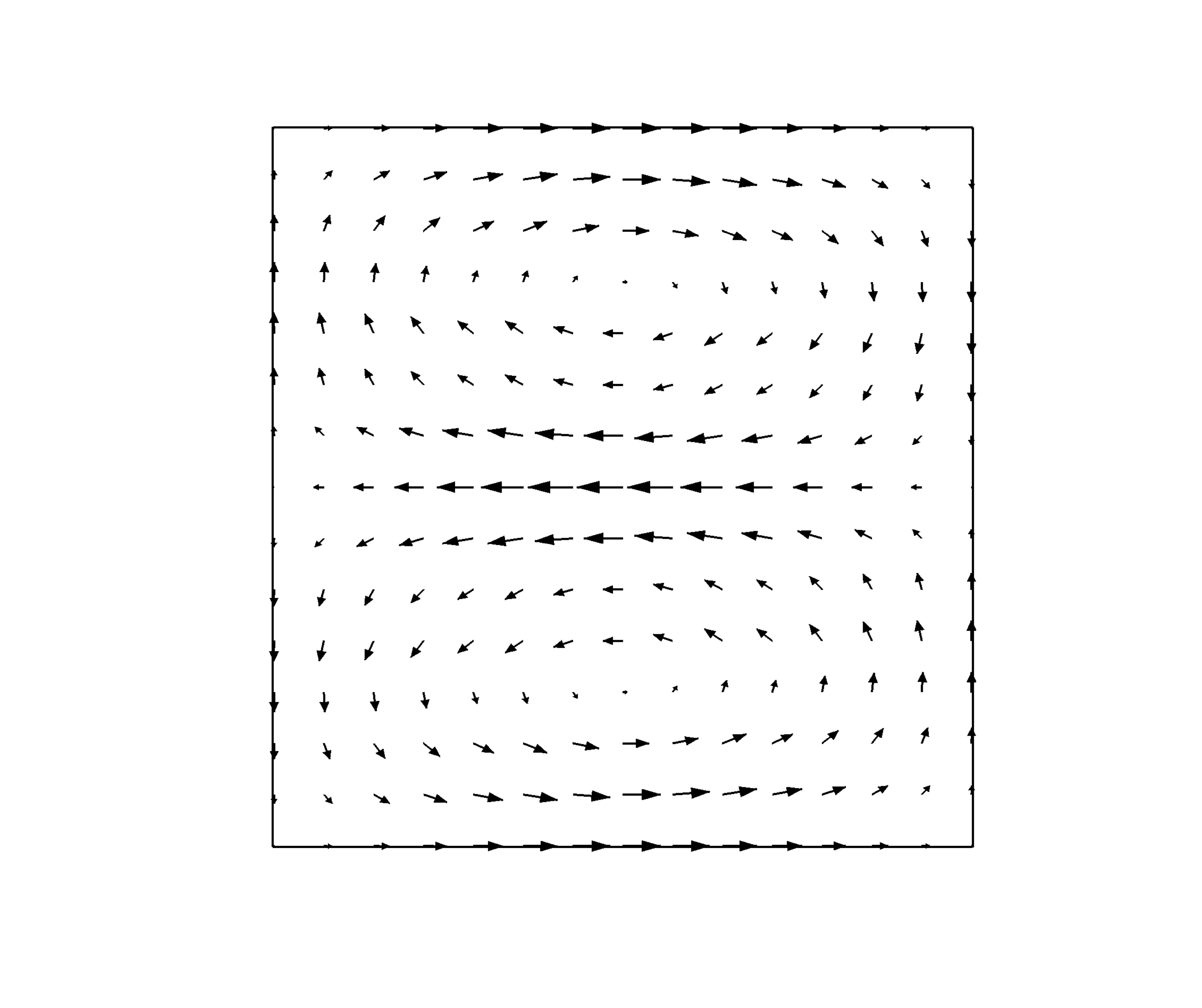}
\hfill
 \includegraphics[trim=18cm 3.58cm 16cm 2.3cm, clip,width=0.49\linewidth]{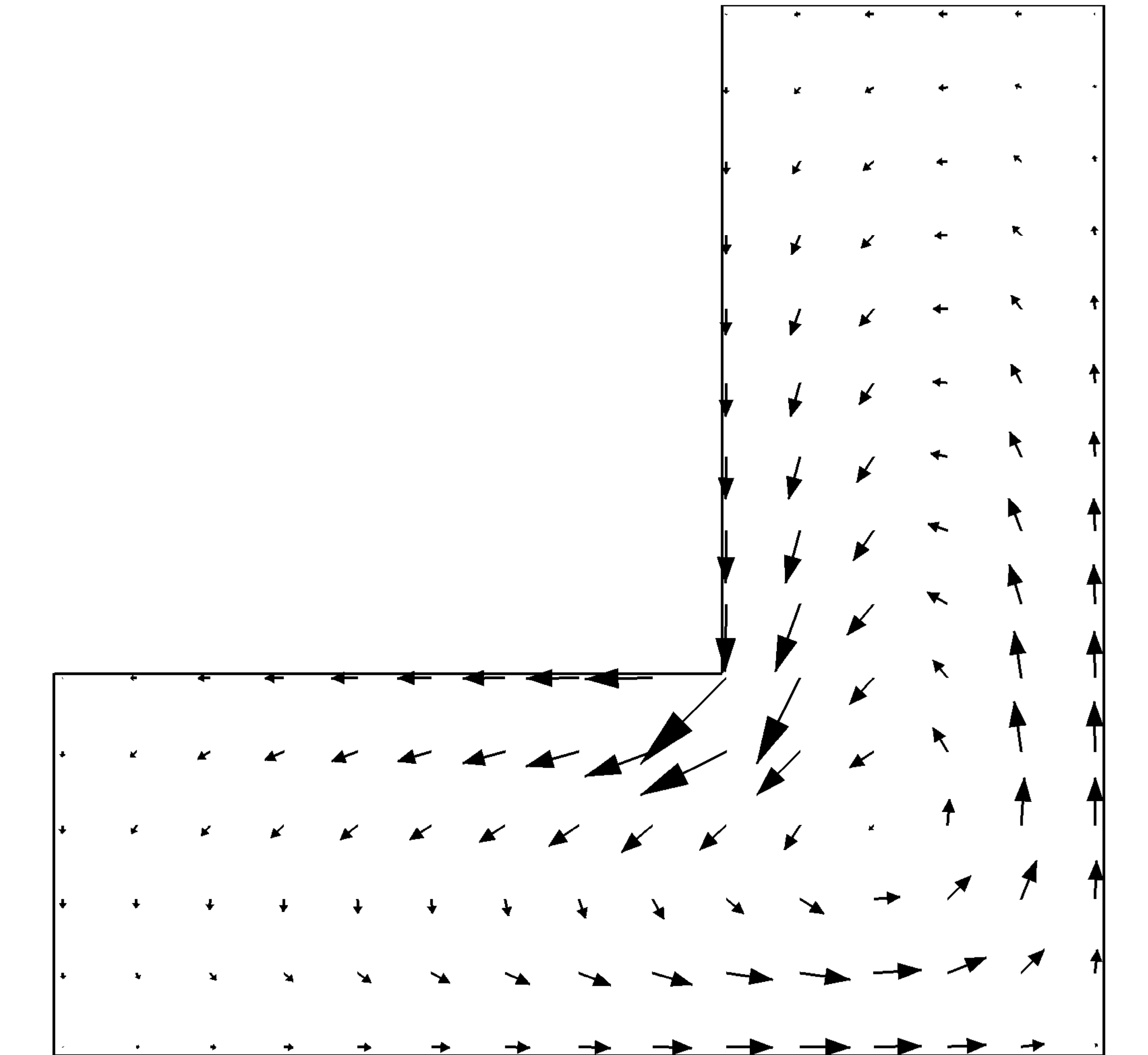}
 \caption{\small Left side: horizontal section of eigenvector $\bff_4$ in the 
 cubic geometry (mode that is not excited). Right side: horizontal section of field $\bsj^3$ for $\lambda/d=4.17$ in the L-shape geometry. \label{modevanish}}
\end{figure}

Now we turn to the description  of the  displacement current $\bsj^3$.
In figure~\ref{fig.sectionJ3}  an horizontal section of this vector field $\bsj^3$ is given for two particular values of $\lambda/d$ (for which $\bsj^3$ is horizontal and very slowly varying in $y_3$).  
On the left, for $\lambda/d=7.5$, a single loop of current is obtained associated with the fundamental eigenvalue $\alpha_1$. 
On the right, for $\lambda/d=3.76$, the second important resonance $\alpha_{24}$ is excited and it appears two concentric loops of current traveled in opposite directions. 
The induced magnetic field $\bsh^3$ does not vanish but has much lower amplitude than in the previous case.

 Eventually in Figure \ref{modevanish}, we draw the horizontal section of eigenvector $\bff_4$ whose strength factor $\int_\Sigma\bfy\wedge\bff_4$ vanishes.
 It consists of two counter rotating loops. No magnetic field is induced and the resonance $\alpha_4$ is not excited from the incident wave.

\paragraph{L-shape geometry.} In that case, the effective permeability tensor $\muef$ is not  scalar anymore. However some symmetries are still  present
and it is possible to show  that   $\muef_{11}=\muef_{22}$ and $\muef_{13}=\muef_{23}=0$. Thus $\muef$ admits  $\{\bfe_1\pm\bfe_2, \bfe_3\}$ as eigenvectors
and its eigenvalues are $\left\{ \muef_{11} \pm \muef_{12}, \muef_{33}\right\}$. 

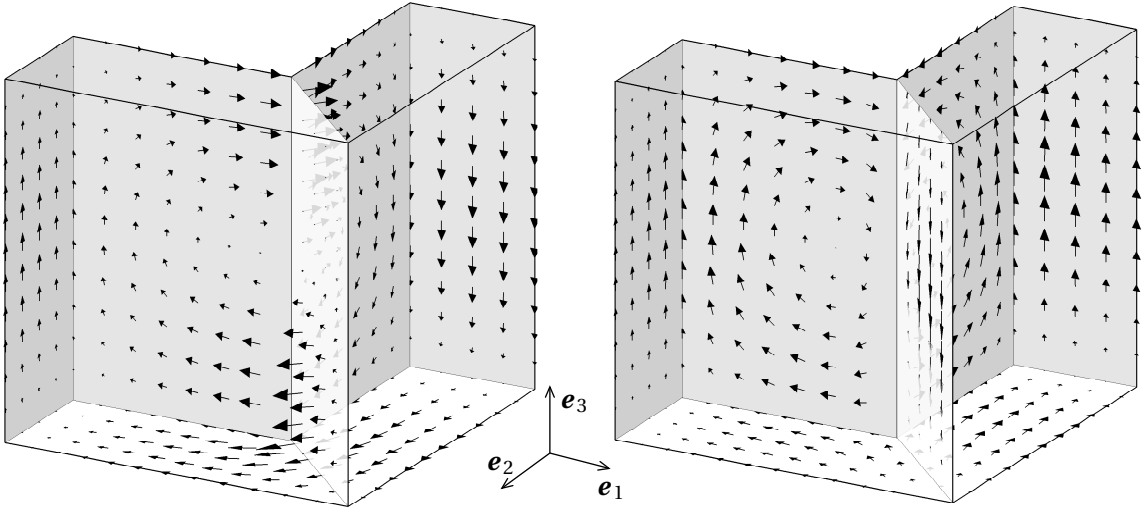
\begin{figure}[!ht]

  \mbox{\hspace{-0.5cm}
\includegraphics[trim=13cm 0cm 13cm 0cm, clip,width=0.5\linewidth]{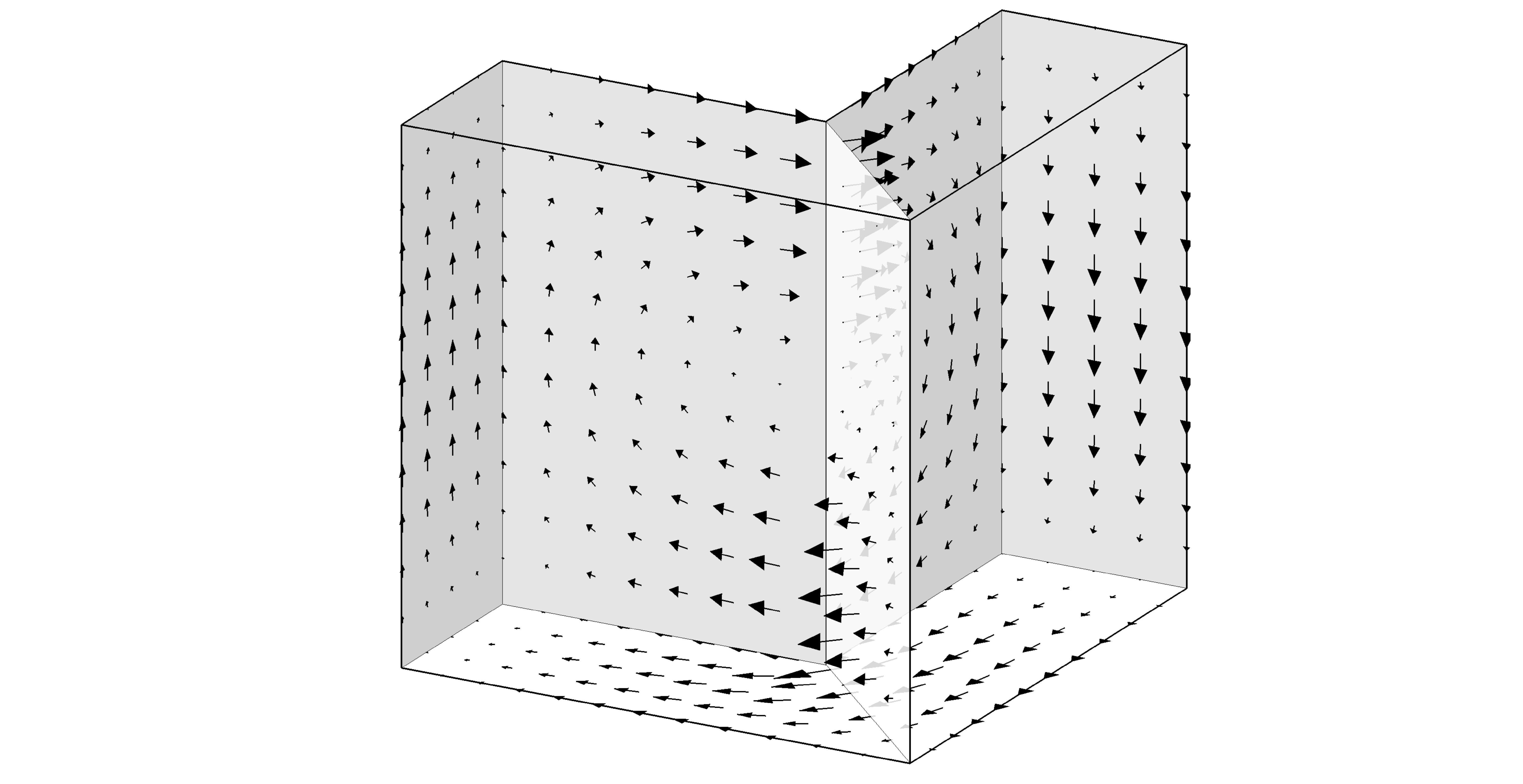}
\hfill
 \includegraphics[trim=13cm 0cm 13cm 0cm, clip,width=0.5\linewidth]{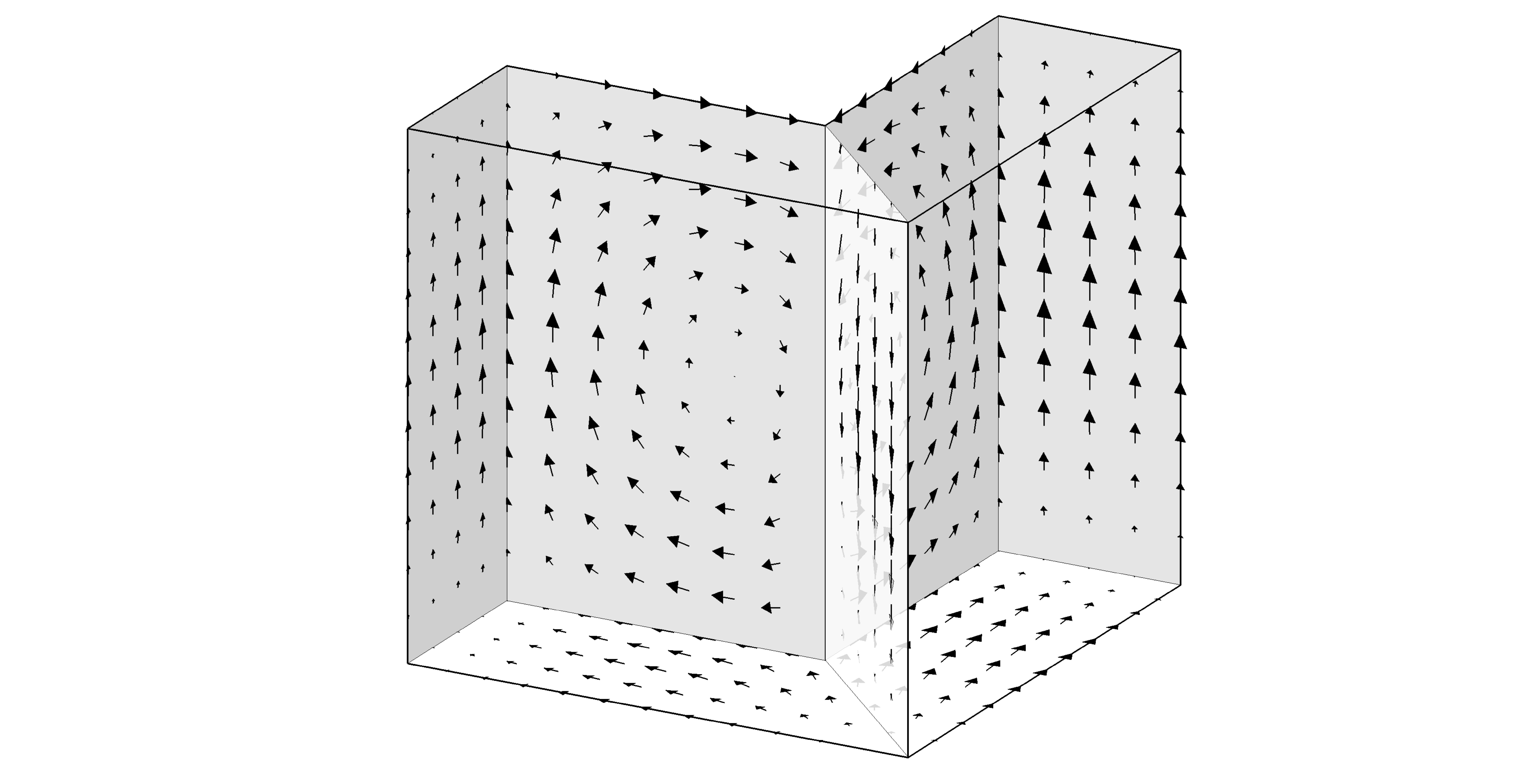}
}
 \begin{textblock*}{2cm}(0.42\linewidth,-18mm)
\begin{tikzpicture}
 \draw [line width=0.5pt, -angle 45] (0,0)-- (0,0.9);
 \draw [line width=0.5pt, -angle 45] (0,0)-- (0.8,-0.21);
 \draw [line width=0.5pt, -angle 45] (0,0)-- (-0.65,-0.48);
\draw (0,0.9) node[anchor=north west] {$\bfe_3$};
\draw (0.8,-0.25) node[anchor=north] {$\bfe_1$};
\draw (-0.65,-0.4) node[anchor=south] {$\bfe_2$};
\end{tikzpicture}
\end{textblock*}
\vspace{-0.2cm}
\caption{\small Representation of field $\bsj^1$ in the case of a L-shape inclusion $\Sigma$. On the left, $\lambda/d=6.2$ (first mode is excited) and on the Right $\lambda/d=5.7$ (second mode is excited) 
 \label{fig.currentJ3}}

\end{figure}

\noindent
In the right hand side  of \ref{fig:mu_1}, the real parts of coefficients $\muef_{11},\muef_{12},\muef_{13}$ are represented 
in term of normalized wavelength $\lambda/d$ varying in the interval $[3,10]$ where  most of the significant resonances are localized.
The two peaks of resonance on the right influence the response to  horizontal  magnetic fields. They are produced by the displacement currents $\bsj^1,\bsj^2$ 
each of them being deduced form the other by a rotation of angle $\pi$ around axis $\bfe_1+\bfe_2$.
In figure \ref{fig.currentJ3} we represent  $\bsj^1$ that we draw only in the back faces of $\Sigma_1$ and in a fictitious surface element located at the junction part of the "L". 
On the left we take  $\lambda/d=6.2$ which is close to the first fundamental mode while on the right $\lambda/d=5.7$ corresponds to the second mode.
In both cases $\bsj^1$ exhibits two loops which  rotate in the same direction in the first case (averaged induced magnetic field $\int_Y\bsh^1$  parallel to $\bfe_1+\bfe_2$ ) and in the opposite direction 
in the second one ($\int_Y\bsh^1$  parallel to $\bfe_1-\bfe_2$).
%
%
%

%
%

 The next significative resonance is obtained for  $\lambda/d$ close to $4.17$ and relies on displacement current  $\bsj^3$.  This vector field
 depicted in the right hand side of figure \ref{modevanish}  is horizontal and independent of $y_3$. The averaged induced magnetic field $\int_Y\bsh^3$ is vertical and contributes  to $\muef_{33}$ only, producing 
a change sign of  $\Re(\muef_{33})$.

\section{Preliminary backgroud}\label{sec:esti_prelim_mu}

\subsection{Two-scale convergence.}

Our study is based on the notion of two-scale convergence which allows to keep a precise description of the electromagnetic field in the periodic microstructure when $\eta$ tends to zero.
We refer to Allaire \cite{allaire} for a precise presentation of the method and we recall here some principal points.

The definitions below are given on a fixed bounded domain $D$ of $\R^3$
 (in most cases $D$ will be a reference ball $B_R$  with radius $R$ so large that $\dom\Subset B_R$).

\begin{definition}[Two-scale convergence] \label{d.twoscale}
We say that a sequence $f_\eta\in L^2(D)$
two-scale converges to $f_0\in L^2(D\times Y)$, and we write $f_\eta(x)\cvd f_0(x,y)$ if, for all $\varphi\in C_c^\infty(D;C^\infty_\sharp(Y))$,
 it holds
\begin{equation}\label{e.stsc}
\lim_{\eta\to 0}\int_{D} f_\eta(x)\, \varphi\xxeta\, \dx =  \int_{D\x Y} f_0(x,y)\, 
\varphi(x,y) \, \dx\,\dy \ .
\end{equation}
The sequence is said  \textit{strongly two-scale convergent} (denoted $f_\eta \cvdf f_0$) if in addition
\begin{equation}\label{two-strong}
\limsup_{\eta\to0}\int_{D} |f_\eta(x)|^2 \dx \ =\ \int_{D\times Y} |f_0(x,y)|^2 \, \dx\, \dy
\end{equation}
\end{definition}
A key justification of definition \eqref{e.stsc} is that any sequence $f_\eta$ which is uniformly bounded in 
$L^2(D)$ admits a two-scale converging subsequence. For such a subsequence, the weak limit exists and is given by the bulk average $\langle f_0\rangle(x) :=\int_Y f_0(x,y)\, dy$.
A consequence of the strong two-scale convergence \eqref{two-strong} is  the following product rule: 
\begin{equation}\label{product}
f_\eta \cvdf f_0 \ ,\ g_\eta \cvd g_0 \quad \Longrightarrow\quad \lim_{\eta\to 0} \int_D f_\eta \, g_\eta \, \varphi \,dx  = \int_D \langle f_0 g_0\rangle\, \varphi\, dx  \ ,
\end{equation}
the latter convergence holding for every $\varphi$ continuous  with compact support. In fact the  assumption on the support of $\varphi$ can be dropped once we know that
 $f_0$ satisfies a suitable  admissibility criterium.  It is the case in particular if $f_0(x,y)=\theta(x)\, \psi(y)$ with $\theta\in L^2(D)$ and $\psi\in L^2_\sharp(Y)$ (see Lemma 5.7 in \cite{allaire}).
 Furthermore the strong convergence $f_\eta \cvdf f_0$ for for  such an admissible  $f_0$ implies that
$$
\lim_{\eta\to0^+}\int_{D}\left| f_\eta(x) - f_0\xxeta\right|^2dx =0.\
$$ 
In particular the strong two-scale convergence $f_\eta \cvdf f_0$ with $f_0(x,y)=f(x)$ implies that $f_\eta \to f$ strongly in $L^2(D)$.

\med
We begin by recalling a classical rule (see \cite{allaire})  which applies to vector valued functions $\bfu_\eta: D \to \C^3$ which are uniformly bounded in $L^2(D)$
as well as  $\ \eta \, \dv \bfu_\eta$. We will use also its straightforward variant in which $\dv \bfu_\eta$ is substituted with $\rt\bfu_\eta$.

\begin{proposition}\label{prop:tsc}
 Let $(u_\eta)_\eta\subset L^2(D;\C^3)$ a sequence of functions such that $u_\eta \cvd \bfu_0$. Then
 \begin{itemize}
 \item[(i)] \ Assume that $\ \eta\, \dv \bfu_\eta \cvd \chi_0\ $. Then for a.e. $x\in D$, it holds\ $\chi_0(x,\cdot) = \dv_y \bfu_0(x,\cdot)$ 
\item[(ii)] \  Assume that  $\ \eta\, \rt \bfu_\eta \cvd \bs{\xi}_0\  $. Then for a.e. $x\in D$, it holds\  $ \bs{\xi}_0(x,\cdot) = \rt_y \bfu_0(x,\cdot)$ 
\end{itemize}
\end{proposition}
Let us stress that the equalities above are meant in the sense of ($Y$- periodic) distributions on $\R^3$.
On the other hand, this result can be localized as follows: let $\Sigma$ be an open subset of $Y$ and $\Sigma_\eta$ defined as in \eqref{Ieta};
then the convergence $\ \eta\, \un_{\Sigma_\eta} \rt \bfu_\eta \cvd  \bs{\xi}_0\  $ implies that $ \bs{\xi}_0(x,\cdot) = \rt_y \bfu_0(x,\cdot)$ in $\mathcal{D}'(\Sigma)$. In the same way we have $ \chi_0(x,\cdot) = \dv_y \bfu_0(x,\cdot)$ in $\mathcal{D}'(\Sigma)$ provided that the convergence
$\ \eta\, \un_{\Sigma_\eta} \dv \bfu_\eta \cvd \chi_0$ holds.

\med
Next we give an improved version of a classical result concerning  a sequence of scalar-valued functions  $u_\eta$ which are uniformly bounded in $W^{1,2}(D)$.
 For every $\eta>0$, let us set 
 $$Y_\eta^k=\eta(Y+k)\quad,\quad I_\eta:=\{ k\in\Z^3\ :\ Y_\eta^k\cap D \not=\emptyset \} .$$
   Then we define the step-wise approximation of $u_\eta$ given by
\begin{equation}\label{def_u_constantmorceau}
[u_\eta]_\eta(x):=\sum_{k\in I_\eta}\Big(\fint_{Y^k_\eta}u_\eta\dx\Big)1_{Y_\eta^k}(x)\ .
\end{equation}
Assume that $u_\eta$  converges weakly to $u$ in $\Hun(D)$. Then by Rellich's Theorem and Poincar\'e inequality, we easily deduce that: 
\begin{equation}\label{u1eta}
 [u_\eta]_\eta\to u\quad \mbox{in}\quad L^2(D)\qquad\text{and}\qquad  v_\eta:= \frac{u_\eta-[u_\eta]_\eta}{\eta} \quad\mbox{is bounded in }L^2(D)\ .
\end{equation}
Possibly passing to a subsequence, we may assume that $v_\eta \cvd v_0$.
 In the following Proposition we make the link  between $v_0$ and the two-scale limit of the sequence $\nabla u_\eta$.
\begin{proposition}\label{example_cvdf}
Let $(u_\eta)$ be a bounded sequence in $W^{1,2}(D)$ such that 
$$  u_\eta \to u   \quad \text{weakly in $W^{1,2}(D)$} \quad ,\quad \frac{u_\eta-[u_\eta]_\eta}{\eta}\cvd v_0\ .$$
Let\ $\ \psi_0(x,y) := v_0(x,y) - \nabla u(x)\cdot (y-[y])\ .$
Then $\psi_0$ belongs to $L^2(D;W^{1,2}_\sharp(Y))$ and it holds
\begin{equation}\label{conv_phi_eta_tilde_bis}
\nabla u_\eta\cvd \nabla_x u(x)+\nabla_y \psi_0(x,y)\ .
\end{equation}
Moreover if $u_\eta$ is independent of $\eta$ (i.e. $u_\eta=u$), then $\psi_0=0$ and we have the strong two-scale convergence $\ \frac{u-[u]_\eta}{\eta} \cvdf \nabla u(x)\cdot (y- [y])$.
\end{proposition}
\begin{remark}{\upshape
 The fact that $\psi_0(x,\cdot)$ belongs to $W^{1,2}_\sharp(Y)$ indicates that the periodic function $v_0(x,\cdot)$ shares  on $\partial Y$ 
 the same jump as the piecewise affine function $\nabla u(x) (y- [y])$.  The existence of $\psi_0$ satisfying 
\eqref {conv_phi_eta_tilde_bis} is proved in \cite{allaire}. Here we give an explicit construction which will be useful in Section 5.} \end{remark}

\med\proof  {\textbf{Step 1}.\ }\ We consider the linear map $A_\eta: W^{1,2}(D) \mapsto L^2(D)$ defined by
$$ A_\eta (u) \ :=\ \frac{u-[u]_\eta}{\eta} - \nabla u(x) \cdot \frac{x-[x]_\eta}{\eta}\ .$$
Since $|x-[x]_\eta| \le \frac{\sqrt{3}}{2} \,\eta $,  it follows from Poincar\'e inequality that for a suitable positive $C$:
\begin{equation}\label{contAeta}
\Vert A_\eta (u) \Vert_{L^2(D)}\ \le C \ \Vert \nabla u \Vert_{L^2(D)}  \qquad  \forall u\in W^{1,2}(D)\ .
\end{equation}
On the other hand, it can be easily checked by using a Taylor expansion that, for $u\in C^2(\ov D)$, $\Vert A_\eta (u) \Vert_{L^2}$ $\le C\, \eta\, .$ 
Therefore by \eqref{contAeta} and by the density of $C^2(\ov D)$ in $W^{1,2}(D)$, we deduce that $A_\eta$ converges strongly to $0$ as $\eta\to 0$.
In particular, for every $u\in W^{1,2}(D)$,  it holds $A_\eta u \to 0$ in $L^2(D)$. Noticing that  $\frac{x-[x]_\eta}{\eta} \cvdf (y-[y])$, it follows 
that $\frac{u-[u]_\eta}{\eta} \cvdf \nabla u(x)\cdot (y- [y])$  (which is the last statement of Proposition \ref{example_cvdf}.


\med 
{\textbf{Step 2}.\ }\ In this step we prove \eqref{conv_phi_eta_tilde_bis} assuming that $u=0$. 
In this case $\psi_0$ coincides with the weak two-scale limit $v_0$ of $\frac{u_\eta-[u_\eta]_\eta}{\eta}$. Possibly passing to a subsequence we may assume that 
\begin{equation}\label{datacvd}
\nabla u_\eta \cvd \xi_0(x,y) \quad ,\quad \frac{u_\eta-[u_\eta]_\eta}{\eta} \cvd v_0(x,y) \ .
\end{equation}
 The two-scale convergence of the whole sequence $\nabla u_\eta$ will be deduced once we can prove that, for a.e. $x\in D$, $\xi_0(x,\cdot )$ agrees with the distributional gradient of $v_0(x,\cdot)$ on $\R^3$. As $v_0(x,\cdot)$ is periodic
this amounts to showing that for a.e. $x\in D$ and for every test function $\bs{\theta}(y)$ in $C^\infty_\sharp(Y;\C^3)$ , it holds 
\begin{equation}\label{nablav0}
- \ \int_Y v_0(x,y)\, \dv_y \bs{\theta} \, dy \ =\  \int_Y \xi_0(x,y) \cdot \bs{\theta}(y)\ dy \ .
\end{equation}
Let $\rho$ an arbitrary localizing function in  $C^\infty_0(D)$. Thanks to an integration  by parts, we obtain
\begin{equation}\label{conv_phi_eta_tilde_bis_preuve}
\int_{D}\nabla u_\eta(x)\rho(x)\cdot\bs{\theta}\xeta\dx=-\int_{D} u_\eta(x)\nabla\rho(x)\cdot\bs{\theta}\xeta\dx-\frac{1}{\eta}\int_{D} u_\eta(x)\rho(x)\dv_y\bs{\theta}\xeta\dx\ .
\end{equation}
 Using the piecewise constant approximation operator $[\cdot]_\eta$ as defined in \eqref{def_u_constantmorceau},  we split the last integral in the right-hand member of \eqref{conv_phi_eta_tilde_bis_preuve} as follows:
\begin{eqnarray}\nonumber
&&\hspace{-0.9cm}\frac{1}{\eta}\int_D u_\eta(x)\rho(x)\dv_y\bs{\theta}\xeta\dx=\int_D\lp\frac{ u_\eta(x)-[ u_\eta]_\eta}{\eta}\rp\rho(x)\dv_y\bs{\theta}\xeta\dx\\\nonumber
&&+\int_D[ u_\eta]_\eta(x)\lp\frac{\rho(x)-[\rho]_\eta}{\eta}\rp\dv\bs{\theta}\xeta\dx+\frac{1}{\eta}\int_D[ u_\eta]_\eta(x)[\rho]_\eta\dv_y\bs{\theta}\xeta\dx\\\nonumber
&&=I_\eta^1+I_\eta^2+I_\eta^3\ .
\end{eqnarray}
We observe that $I_\eta^3$ vanishes : indeed, the functions $[ u_\eta]_\eta$ and $[\rho]_\eta$ are constant on each cell $Y_\eta^k$ where  $\dv_y\bs{\theta}$ has a vanishing mean value. 
On the other hand, since $[ u_\eta]_\eta$ strongly converges to $u=0$ in $L^2(D)$, we infer that $I_\eta^2\to 0$ (indeed $\frac{\rho-[\rho]_\eta}{\eta}$ is bounded in $L^2(D)$ whereas $\dv\bs{\theta}\in L^\infty(D)$).

Eventually passing to the limit $\eta\to 0$ in \eqref{conv_phi_eta_tilde_bis_preuve} and taking into account the two-scale convergences \eqref{datacvd} , we obtain the equality
\begin{equation*}
 \int_{D\x Y} \xi_0(x,y)\, \rho(x)\cdot\bs{\theta}(y)\, dxdy\ =\ -\int_{D\x Y}v_0(x,y)\, \dv_y\, \bs{\theta}(y)\,  \rho(x)\, dxdy\ ,
\end{equation*}
from which follows the relation \eqref{nablav0} by the arbitrariness of the test function $\rho$. 

\med 
{\textbf{Step 3}.\ }\ We consider now  a general sequence $u_\eta$ converging weakly to $u$ and we apply Step 2 to the translated sequence $\tilde u_\eta= u_\eta-u$.
Then  $\tilde v_\eta:=\frac{\tilde u_\eta-[\tilde u_\eta]_\eta}{\eta}= v_\eta - \frac{ u-[ u]_\eta}{\eta}$.  Thus, by Step 1, the two-scale limit $\tilde v_0$ of $\tilde v_\eta$  is given by
$\tilde v_0(x,y)  = v_0(x,y)-\nabla u(x) \cdot(\bfy-[\bfy]) = \psi_0(x,y)$. 
On the other hand, as $\tilde u_\eta\to 0$, we know by Step 2 that $\nabla\tilde u_\eta\cvd \nabla_y\tilde v_0$. It follows that $$\nabla u_\eta=\nabla\tilde u_\eta+\nabla u\cvd \nabla_y \psi_0(x,y) +\nabla u(x)\, .$$
The proof of Proposition \ref{example_cvdf} is finished. 
\qed

%

%
%
%

\subsection{Geometric averaging.} \label{sec:circ}

Since the advent of gauge theories, it is now a classical idea in Physics to see magnetic fields as differential 1-forms in $\R^3$. In our case it turns out that
by Lemma \ref{p.magnetic} the magnetic field $\bsh_0(x,\cdot)$ is curl-free in the simply-connected domain $\Sigma^*$ (closed 1-form). 
It follows that the circulation of $\bsh_0(x,\cdot)$ along curves in $\Sigma^*$ joining opposite points of $\partial Y$ is independent of the line 
and of the end points we chose.
It is then natural to define the ``mean circulation'' vector, denoted by $\oint\bsh_0$, and characterized by:
\begin{equation}\label{def_circ}
\lp \oint\bsh_0(x,\cdot) \rp\cdot\bfe_k:=\int_{\Y_k}\bsh_0(x,\cdot) \cdot\bfe_k\,d\mathcal{H}^1\ ,
\end{equation}
where $\Y_k\subset\Sigma^*$ is a curve joining two opposite points on the faces of $\partial Y$ orthogonal to $\bfe_k$.

\textit{A priori} this definition makes sense only for \textit{regular} functions. It can be extended to functions in $\Hunper(Y)$ as follows:
\begin{lemma}\label{lem:def_circulation}
Let  $\bfu\in L^2_\sharp(Y;\C^3)$ such that\ $\rt\bfu=0$ in $\Sigma^*$. 
Then there exists a  \textit{unique} vector $\oint\bfu\in\C^3$ and a function $\psi\in \Hunper(\Sigma^*)$ (unique up to a constant) such that
\begin{equation}\label{def:circ2}
\bfu=\nabla\psi+\oint\bfu\quad\mbox{in}\quad \Sigma^*\ .
\end{equation}
Moreover this circulation vector satisfies the following properties:
\begin{itemize}
\item[i)] For all $\bff\in L^2_{\sharp}(Y;\C^3)$ such that $\dv\bff=0$ in $Y$ and $\bff=0$ in $\Sigma$ we have
\begin{equation}\label{def:circ}
  \int_Y\bfu\cdot\bff\dy=\oint\bfu\,\cdot\int_Y\bff\dy\ .
\end{equation}

\item[ii)] If $\bfu$ is continuous, then for $\bfk\in\{1,2,3\}$
\begin{equation}\label{def:circ3}
\lp\oint\bfu\rp\cdot\bfe_k=\int_0^1\bfu\big(\gamma(s)\big)\cdot\gamma'(s)\, ds\ ,
\end{equation}
for all $\gamma\in C^1([0,1];\Sigma^*)$ such that $\gamma^k(1)-\gamma^k(0)=\bfe_k$.
\item[iii)] There exists a constant $C>0$ such that
\begin{equation}\label{esti_circ}
 \Big|\oint\bfu\Big|\leq C\|\bfu\|_{L^2(Y;\C^3)}\ .
\end{equation}
\end{itemize}
\end{lemma}
\begin{remark}\label{rem_circ}\nshape{Note that if  $\rt\bfu=0$ in whole $Y$, then\, 
$\oint \bfu=\int_Y\bfu\,  $ and we recover from \eqref{def:circ} a classical property (linked to "div-curl Lemma").
 On the other hand, by applying \eqref{def:circ} to $\bff:= \bfv \wedge \bfz$ with $\bfz$ varying over $\C^3$, we get the following  variant of \eqref{def:circ}: 
\begin{equation}\label{wedge}
  \int(\bfu\wedge \bfv) \dy=\oint\bfu\,\wedge \int\bfv\dy\ , \qquad \text{if \quad $\rt\bfv=0$ \ ,\  $\bfv=0$ in $\Sigma$}\ .
  \end{equation} 

} \end{remark}
\proof 
As $\Sigma^*$ is simply-connected, the existence of the unique vector $\oint\bfu\in\C^3$ and of a function $\psi \in\Hunper(\Sigma^*)$ such that \eqref{def:circ2} holds 
in the case of a smooth function $\bfu$ is a classical result in differential geometry (see for instance \cite{berger_gaostiau} p. 197). 
The extension to functions $\bfu\in L^2(Y,\C^3)$ such that $\rt \bfu=0$ in $\Sigma^*$ follows by using a density argument. 

\med
Now we fix  $\bff\in L^2_{\sharp}(Y;\C^3)$ such that $\dv\bff=0$ in $Y$ and $\bff=0$ in $\Sigma$. Moreover we consider 
 $\tilde\psi\in\Hunper(Y)$ to be the harmonic extension of $\psi$ in $\Sigma$. Since $\bff=0$ in $\Sigma$, we have
 $$ 
 \int_Y\bfu\cdot\bff\dy\ = \int_{\Sigma^*}(\nabla \psi + \oint\bfu)\cdot\bff\dy\ =\  \int_Y \nabla \tilde\psi \cdot \bff+ \oint\bfu\cdot\int_Y\bff\dy= \oint\bfu \cdot \int_Y \bff\, dy\ ,
$$
where the last equality holds since $\dv\bff=0$. Relation \eqref{def:circ} is proved.

\med
For $\bsu$ being continuous, the associated potential $\psi$ is Lipschitzversions dŽfinitives on
$\Sigma^*$ and we get
$$
 \int_0^1\bfu\big(\gamma(s)\big)\cdot \gamma'(s)\, ds
= \int_0^1 (\nabla \psi+ \oint\bfu)\big(\gamma(s)\big)\cdot \gamma'(s)\, ds =
\psi \big(\gamma(1)\big)-\psi\big(\gamma(0)\big) + \oint\bfu\cdot \bfe_k \ .
$$
The equality \eqref{def:circ3} follows thanks to the periodicity of $\psi$.

\med 
For the proof of  \eqref{esti_circ} we consider three cylinders $C_k\subset Y$ for $k\in\{1,2,3\}$.
Each of them is in direction $\bfe_k$, cross entirely the unit cell $Y$ and is such that $\Sigma\cap C_k=\emptyset$.
We introduce functions $\bff_k\in L_\sharp^2(Y;\C^3)$ given by $\bff_k(y)=\frac{1}{|C_k|}1_{C_k}(y)\bfe_k$. By construction 
those functions are admissible in \eqref{def:circ} and satisfy $\int_Y\bff_k\,\dy=\bfe_k$. Thanks to \eqref{def:circ} and to the Cauchy-Schwartz inequality it comes
$$
\Big|\oint\bfu\cdot\bfe_k\Big|=\Big|\oint\bfu\cdot\int_Y\bff_k\Big|=\Big|\int_Y\bfu\cdot\bff_k\Big|\leq \frac{1}{\sqrt{|C_k|}}\|\bfu\|_{L^2(Y;\C^3)}\ ,\qquad \mbox{for }k\in\{1,2,3\}.
$$
The estimate \eqref{esti_circ} follows.
\qed

\subsection{Miscellaneous results.}
  For the commodity of the reader we recall a classical result related Sobolev space $\Hunper(Y;\C^3)$ 

\begin{lemma}\label{prop:equiv_norm}
Let  $\bfu\in L^2_\sharp(Y;\C^3)$ such that \ $\rt\bfu\in L^2_\sharp(Y;\C^3)$ \  and \  $\dv\bfu\in  L^2_\sharp(Y)$. Then $\bfu\in \Hunper(Y;\C^3)$ and
\begin{equation}\label{equiv_norm_per}
\int_Y |\nabla \bfu|^2 \ = \ \int_Y | \rt \bfu|^2  + \int_Y |\dv \bfu|^2 \ .
\end{equation}
\end{lemma}
\proof
As elements of $L^2_\sharp(Y)$, the functions $\bfu$, $\rt\bfu$ and $\dv\bfu$ can be developed in Fourier series. There exists therefore $\la\bfc_k\in\C^3\ :\ \bfk\in\Z^3\ra$ such that
$$
\bfu(\bfy)=\sum_{\bfk\in\Z^3}e^{2i\pi\bfk\cdot\bfy}\bfc^k\ ,\qquad
\rt\bfu(\bfy)=\sum_{\bfk\in\Z^3}2i\pi e^{2i\pi\bfk\cdot\bfy}\ \bfk\wedge \bfc^k\ ,\qquad
\dv\bfu(\bfy)=\sum_{\bfk\in\Z^3}2i\pi e^{2i\pi\bfk\cdot\bfy}\ \bfk\cdot\bfc^k\ .
$$
with 
$$ 
\sum_{\bfk\in\Z^3}|\bfc^k|^2=\|\bfu\|_{L^2(Y;\C^3)}^2\ , \quad \sum_{\bfk\in\Z^3}|\bfk\wedge\bfc^k|^2= 4 \pi^2\int_Y | \rt \bfu|^2\ ,\qquad\sum_{\bfk\in\Z^3}|\bfk\cdot\bfc^k|^2= 4 \pi^2\int_Y | \dv \bfu|^2\ .
$$
Since $|\bfk\wedge\bfc^k|^2+|\bfk\cdot\bfc^k|^2 = |\bfk|^2 |\bfc_k|^2$, we infer that $\sum_{\bfk\in\Z^3} |\bfk|^2 |\bfc_k|^2 <+\infty$, thus $\bfu$ belongs to $\Hunper(Y;\C^3)$ 
Moreover relation \eqref{equiv_norm_per} follows obviously.

\qed


\med
In order to establish the strong convergence in $L^2_{\mathrm{loc}}$ of the field $(\bse_\eta,\bsh_\eta)$, we will use the celebrated div-curl Lemma below for which we refer to
\cites{Murat,tartar}.
\begin{lemma}\label{lem:divrot} Let $B$ be a bounded domain of $\R^3$. Let $(\bfu_\eta)$ and $(\bfv_\eta)$ two sequences of $L^2(B;\R^3)$ converging respectively to $\bfu$ and $\bfv$ weakly in $L^2(B;\R^3)$. If moreover we have 
$$
(\dv \bfu_\eta)_\eta\quad\mbox{is compact in }W^{-1,2}(B)\qquad\text{and}\qquad
(\rt \bfv_\eta)_\eta\quad\mbox{is compact in }(W^{-1,2}(B))^3\ ,
$$
then the sequence $(\bfu_\eta\cdot\bfv_\eta)_\eta$ converges to $\bfu\cdot\bfv$ in the distributional sense in~$B$.
\end{lemma}
For checking the strong compactness assumption in $W^{-1,2}(B)$, the following characterization will be useful:
\begin{lemma}\label{prop:comp_fort_Hmoins1}
  Let $(T_n)_n\subset W^{-1,2}(B)$ a bounded sequence. Then $(T_n)$ is relatively compact (for the topology of the norm) if and only if, for all sequence $(\varphi_n)_n $  such that $\varphi_n\cvw0$ in  $\Hun_0(B)$, we have
$\lim_{n\to+\infty}<T_n,\varphi_n>=0\ .$
\end{lemma}
%

\section{Two-scale analysis of the electromagnetic  field }\label{sec:electromagnetic}

In this section we fix a ball  $B_R$ such that $\dom\Subset B_R$. We  are going to identify the two-scale limit of the electromagnetic field $(\bse_\eta,\bsh_\eta)$
under the assumption that it is uniformly bounded in $L^2(B_R)$. In fact this analysis includes the divergence-free vector field $\bsj_\eta$ defined below
which represents a normalized version of the so-called "displacement current":
\begin{equation}\label{Jeta}
\bsj_\eta(x)=\eta\eps_\eta\bse_\eta
\end{equation}
Along this section, we will make the following hypothesis
\begin{equation}\label{born_prio1}
 \sup_{\eta>0}\Big(\|\bse_\eta\|_{L^2(B_R)}+ \|\bsh_\eta\|_{L^2(B_R)}+\|\bsj_\eta\|_{L^2(B_R)}\Big)<+\infty\ .
\end{equation}
Thanks to the estimate \eqref{born_prio1}  which will be established \textit{a posteriori} in Section \ref{sec:proof_prin},
we may assume, possibly after extracting subsequences, that it holds
\begin{equation}\label{twoscalconv}
 \bse_\eta\cvd \bse_0\ ,\qquad \bsh_\eta\cvd \bsh_0\quad\mbox{and}\quad \bsj_\eta\cvd \bsj_0\ .
\end{equation}
for suitable  $\bse_0$, $\bsh_0$ and $\bsj_0$ belonging to $L^2(B_R\x Y;\C^3).$

\med
In view of the convergences in \eqref{twoscalconv}, we define the effective electromagnetic field $(\bse,\bsh)$ to be
\begin{equation}\label{effEH}
\bse(x) \ := \int_Y \bse_0(x,y)\, dy \quad,\quad \bsh(x) \ := \oint \bsh_0(x,\cdot)\ .
\end{equation}
Recall that under \eqref{twoscalconv}, $\bse(x)$ represents the weak limit of  $\bse_\eta$ (bulk average)
whereas  $\bsh(x)$ is associated with the new averaging procedure introduced in Subsection \ref{sec:circ}.

\med
In Subsection \ref{sec:electric}, we identify $\bse_0(x,\cdot)$  for $x\in \dom$.
The characterization of vector fields $\bsh_0(x,\cdot), \bsj_0(x,\cdot)$ for $x\in \dom$ is a quite involved issue developed along Subsection \ref{sec:magnetic}. 
Besides, for $x\notin\dom$, $\bsj_0(x,\cdot)$ vanishes and  it is a direct consequence of the uniform convergence issue developed in Subsection \ref{subsec:far} that
 the fast oscillations of $(\bse_\eta,\bsh_\eta)$ disappear at a positive  distance from $\partial\dom$.
 Accordingly, $\bse_0(x,\cdot), \bsh_0(x,\cdot)$ are constant and agree with their respective averages:
  \begin{equation}\label{EHconstant}
 \bse_0(x,\cdot)  = \bse(x) \quad , \quad \bsh_0(x,\cdot) = \bsh(x)  \qquad \text{for \, a.e. $x\in B_R\setminus\dom$}\ .
\end{equation}

\subsection{Oscillating electric field and effective permittivity.}\label{sec:electric}
Our aim here is to  identify, for $x\in\dom$, the two-scale limit $\bse_0$ given in \eqref{twoscalconv} in term of its bulk average $\bse(x)$ (see \eqref{effEH}).
First we show that $\bse_0(x,\cdot)$ solves an electrostatic problem in the unit cell:
\begin{lemma} 
 For almost every $x\in \dom$, the periodic function $\bfu:=\bse_0(x,\cdot)$ satisfies in the distributional sense in $Y$
\begin{equation}\label{PE0}
\rt_y \bfu=0\ ,\qquad
\dv_y \bfu=0\ \ \mbox{in }\  \dom\x  \Sigma^*\,,\quad\mbox{and}\quad\bfu=0\ \ \mbox{in }\  \dom\x \Sigma\ .
\end{equation}
In particular the restriction of $\bse_0(x,\cdot)$ to $\Sigma^*$ belongs to $W^{1,2}_\sharp(\Sigma^*;\C^3).$ 
\end{lemma}
\proof 
By \eqref{born_prio1} and the first  equation in \eqref{prb:prin_mu}, we infer that $ \eta\, \rt \bse_\eta \to 0$ in $L^2(\dom)$. Thus by applying
Proposition \ref{prop:tsc} (ii), we deduce that\  $\rt_y \bse_0(x,\cdot)$ vanishes everywhere.

\med
From the second Maxwell equation in \eqref{prb:prin_mu}, we deduce that $\dv(\eps_\eta\bse_\eta)=0$ in $\dom$. Since $\eps_\eta$ is constant in 
$\dom\setminus \Sigma_\eta$, it follows in particular that $\eta\dv(\bse_\eta)=0$ in $\dom\setminus\Sigma_\eta$. 
Passing to the limit $\eta\to 0$ with the help of assertion (i) of Proposition \ref{prop:tsc} (and of the localized version indicated after the statement), we derive the equality  $\dv_y \bse_0=0$ holding in $\mathcal{D}'(Y\setminus\Sigma)$.

\med
Eventually, we observe that  the equality $\eta\, \bsj_\eta= \eps_r \bse_\eta$  holds in $\Sigma_\eta$. By the bound \eqref{born_prio1}, we deduce that $\bse_\eta \un_{\Sigma_\eta} \to 0$ in $L^2(\dom)$.
Accordingly  $\bse_0(x,\cdot)$ vanishes in $\Sigma$. 
\qed

\med
We are now able to characterize the set of solutions of the cell problem \eqref{PE0}.
\begin{proposition}\label{prop:decomp_E0}
The set of solutions of \eqref{PE0} is a three dimensional vector space spanned by the fields $\bse^k\in L^2(Y;\R^3)$ defined for $k=1,2,3$ by :
\begin{equation}\label{def:Ei}
 \bse^k=\bfe_k+\nabla_y\chi_k\ .
\end{equation}
Functions $\chi_k$  are the unique elements of \  $W^{1,2}_\sharp(Y)$  satisfying \eqref{localE}, that is:
\begin{equation*}
\Delta_y\chi_k=0\quad\mbox{in}\quad \Sigma^*\qquad\mbox{and}\qquad \chi_k=-y_k\quad\mbox{in}\quad \Sigma\ .
\end{equation*}
In particular, $\bse_0$ can be decomposed as
\begin{equation}\label{decomp_E0}
\bse_0(x,y)=\sum_{k=1}^3 E_k(x)\, \bse^k(y)\ ,
\end{equation}
being  $\ \disp E_k(x) = \int \bse_0(x,y)\cdot\bfe_k\, dy$.

\end{proposition}
\proof
By construction $\bse_0$  given  in  \eqref{decomp_E0} satisfies \eqref{PE0} 
and the average condition $\int_Y\bse_0(x,y)\dx=\bse(x)$. 
More generally, if $V$ denotes the subspace of $\Hunper(Y)$ consisting of all solutions to \eqref{PE0},  we see that the linear 
 map $ \bfu \in V \mapsto \int_Y \bfu \, dy \in \C^3$
is surjective. The uniqueness as well as the fact that $\mathrm{ dim}(V) =3$ will follow if we show that this map is injective.
\med
Let $\bfu \in V$ such that $\int_Y \bfu=0$.  Then as $\rt \bfu=0$ in $Y$, there exist $\psi\in W^{1,2}_\sharp(Y)$
such that $\bfu =\nabla\psi$ in $Y$. Thanks to  \eqref{PE0}, it holds $\Delta\psi=0$ in $\Sigma^*$ while $\nabla \psi=0$ in $\Sigma$. As $\Sigma$ is connected, 
we may assume that $\psi$ vanishes in $\Sigma$ as well as its trace on $\partial\Sigma$. It follows that $\psi$ vanishes everywhere, thus  $\bfu=0$.
\qed
%

\paragraph{Effective permittivity tensor.}
The functions $\bse^k(y)$ defined in \eqref{def:Ei} that we call \textit{shape electric fields}  depend only on the geometry of $\Sigma$. They determine the oscillating
behavior of the electric field as well as the effective permittivity tensor $\epsef$ given in \eqref{def_epsef_intro} (up to 
the positive real permittivity factor $\eps_e$ which represents the permittivity of the matrix). Indeed the tensor given in \eqref{def_epsef_intro} can be written as follows:
\begin{equation}\label{epsef2}
\epsef_{kl}:=\eps_e\int_Y\bse^k\cdot\bse^l\ ,
\end{equation}
In particular we see that this tensor is real symmetric, positive-definite and independent on the frequency.
More precisely, by using Jensen's inequality, we find that
\begin{equation}\label{positive-eps}
\epsef z\cdot z \ =\ \eps_e \, \int_Y | \sum_{k=1}^3 z_k\, \bse^k |^2 \ \ge \ \eps_e \, |z|^2 \qquad \text{for all $z\in \C^3$}
\end{equation}
As a consequence no specific effect (resonances, dispersion) is expected concerning the permittivity law of the composite structure under study.

\subsection{Oscillating magnetic field and displacement current.}\label{sec:magnetic}

In this subsection, we are going to identify  the two-scale limits $\bsh_0$ and $\bsj_0$ given in \eqref{twoscalconv} in terms of the effective magnetic field $\bsh(x)$
defined in \eqref{effEH}. We start with the following result:
\begin{lemma}\label{p.magnetic}
For almost all $x\in \dom$, the periodic vector fields $\bsh_0(x,\cdot)$ and $\bsj_0(x,\cdot)$ satisfy:
\begin{equation}\label{PH0}
 \dv_y \bsh_0(x,\cdot) \ =\ 0  \quad ,\quad \rt_y \bsh_0(x,\cdot) \ =\ -i\, \omega\, \eps_0\, \bsj_0(x,\cdot)  \quad \mbox{\upshape{ in $\mathcal{D}'(\R^3)$}}  \ ,\end{equation}
 \begin{equation}\label{PJ0}
  \bsj_0(x,\cdot) \ =\ 0  \quad \mbox{\upshape{a.e. in $\Sigma$}} \quad ,\quad \rt_y \bsj_0(x,\cdot) \ =\ i\omega\mu_0 \eps_r\, \bsh_0(x,\cdot)  \quad \mbox{\upshape{ in $\mathcal{D}'(\Sigma)$}} \ .
  \end{equation}
 In particular $\bsh_0(x,\cdot)$ is an element of $W^{1,2}_\sharp(Y,\C^3)$.
 
\end{lemma}
%
%
\proof Recalling that $\bsh_\eta$ is divergence-free, the first relation in \eqref{PH0} follows from the assertion (i) of Proposition \ref{prop:tsc}.
Next by the second equation in \eqref{prb:prin_mu} and \eqref{Jeta}, we have 
$$
\eta\rt\bsh_\eta=-i\omega\eps_0\eta\eps_\eta\bse_\eta=-i\omega\eps_0\bsj_\eta\ .
$$
By taking the two-scale limit of the left-hand member exploiting the assertion (ii) of Proposition \ref{prop:tsc}, we  are led to the second equation in \eqref{PH0}.
In particular we find that  $\rt \bsh_0(x,\cdot), \dv \bsh_0(x,\cdot)$ are elements $L^2_\sharp(Y)$.  Therefore, by Lemma \ref{prop:equiv_norm}, we  conclude that 
$\bsh_0(x,\cdot)$ belongs to $W^{1,2}_\sharp(Y,\C^3)$.

\med Next we observe that, by construction, the vector field $\bsj_\eta(x)$ vanishes outside $\Sigma_\eta$, thus obviously  $\bsj_0(x,\cdot) =0$ a.e. in $\Sigma^*$.
 On the other hand, by \eqref{def_epseta} and \eqref{Jeta}, we have $ \eta \, \bsj_\eta=  \eps_r\, \bse_\eta$ in $\Sigma_\eta$ so that by the first equation in \eqref{prb:prin_mu},   it holds
 $$ \eta \rt \bsj_\eta\  =\  i\,\omega\, \mu_0 \,\eps_r \, \bsh_\eta   \quad \text{ in $\Sigma_\eta$}\ .$$
Passing to the limit $\eta\to 0$ with the help of assertion (ii) of Proposition \ref{prop:tsc} (see the localization argument after the statement), we derive
the equality  $\rt_y \bsj_0(x,\cdot) \ =\ i\omega\mu_0 \eps_r\, \bsh_0(x,\cdot)$ holding in  $\mathcal{D}'(\Sigma)$.
The proof of Lemma \ref{p.magnetic} is complete. \qed

\med 
We notice that, by eliminating $\bsj_0$ in relations \eqref{PH0} \eqref{PJ0},  we obtain for $\bsh_0$ a system of equations similar to that obtained in \eqref{PE0} for $\bse_0(x,\cdot)$:
\begin{equation}\label{PH02}
\dv_y\bsh_0=0 \ ,\quad \text{in $Y$}\quad\,\quad \rt_y \bsh_0=0\quad\text{in $\Sigma^*$}\quad,\quad  \Delta_y\bsh_0+\eps_r\,k_0^2\,\bsh_0=0\quad\text{in $\Sigma$}
\end{equation}
However, as we have no information on the possible tangential jump of $\bsj_0(x,\cdot)$ across $\partial\Sigma$,  it is non straightforward in this case to see  that the set of $W^{1,2}_\sharp(Y,\C^3)$
of solutions to \eqref{PH02} is
still of dimension three or equivalently  that the solution to the system is unique for a given average $\bsh(x)$ (see \eqref{effEH}).
This apparent difficulty is  overcome once we show that $\bsh_0(x,\cdot)$ satisfies  a suitable variational principle.
This is done in the next crucial Lemma.

\begin{lemma} \label{variational}
Let $\bfv \in W^{1,2}_\sharp(Y,\C^3)$ such that $\rt \bfv =0$ in $\Sigma^*$ and \, $\disp\oint \bfv=0$. Then for almost all $x\in \Omega$, the
periodic field $\bfw=\bsh_0(x,\cdot)$ satisfies the equation
\begin{equation} \label{weakform}
\int_\Sigma \rt \bfw\cdot\rt \bfv-\eps_r\,k_0^2\,\int_Y\bfw\cdot\bfv=0\ .
\end{equation}
\end{lemma}
\proof
From the Maxwell system \eqref{prb:prin_mu} and taking into account the definition of $\eps_\eta$  in \eqref{def_epseta} and  $k_0^2=\eps_0\,\mu_0\,\omega^2$,
it is straightforward to deduce that the magnetic field $\bsh_\eta$ satisfies the variational equality
\begin{equation}\label{variationalH}
\frac{1}{\eps_e} \int_{\dom\setminus \Sigma_\eta}\rt\bsh_\eta\cdot\rt\bphi_\eta + \frac1{\eps_r} \int_{\Sigma_\eta}\eta^2 \rt\bsh_\eta\cdot\rt\bphi_\eta\ =\ k_0^2\int_{\dom}\bsh_\eta\cdot\bphi_\eta\ .
\end{equation}
holding for every smooth vector field $\bphi_\eta$ compactly supported in $\dom$. 

\med
Let $\rho(x)\in C^\infty_c(\dom)$ and plug\,  $\bphi_\eta(x):=\rho(x)\bfv(x/\eta)$ as a test function in \eqref{variationalH} with $\bfv$  satisfying the assumptions of the Lemma.
Then, as $\bfv$ is curl-free in $\Sigma^*$, it holds  
$$ \rt(\bphi_\eta)=  \nabla\rho(x)\wedge\bfv\xeta  + \frac{\rho(x)}{\eta}\, \left(\rt_y \bfv\right) \xeta\ .$$
It follows that we have the following strong two-scale convergences 
$$  \bphi_\eta  \cvdf   \rho(x)\, \bfv(y) \quad,\quad \eta\, \rt(\bphi_\eta)\cvdf   \rho(x)\, \rt_y \bfv  \ .$$
Then, by applying \eqref{product} and recalling that $\bsh_\eta \cvd\bsh_0$ and $\eta \rt \bsh_\eta \cvd \rt_y \bsh_0$ (see Proposition \ref{prop:tsc}),  we have
$$
\lim_{\eta\to 0}  \int_{\Sigma_\eta}\eta^2 \rt\bsh_\eta\cdot\rt\bphi_\eta = \int_{\dom\times\Sigma} \rho(x)\, \rt_y\bsh_0\cdot\rt_y\bfv \quad,
\quad \lim_{\eta\to 0} \int_{\dom}\bsh_\eta\cdot\bphi_\eta
= \int_{\dom\times\Sigma}\rho(x)\, \bsh_0\cdot\bfv\ .
$$
In view of \eqref{variationalH}, it is then enough to prove the following claim 
\begin{equation}\label{claim19}
\lim_{\eta\to 0} \int_{\dom\setminus \Sigma_\eta}\rt\bsh_\eta\cdot\rt\bphi_\eta \ =\ 0\ .
\end{equation}
 Indeed thanks to the convergences above, by passing to the limit in \eqref{variationalH} , we infer that
 $$\int_{\dom\x\Sigma}\rho(x)\,\rt_y\bsh_0(x,y)\cdot\rt\bfv(y)\, dxdy
=k_0^2 \eps_r \int_{\dom\x Y}\rho(x)\, \bsh_0(x,y)\cdot\bfv(y)\, dxdy\  ,$$
and \eqref{weakform} follows by the arbitrariness of the test function $\rho(x)$.

 \med
 To prove claim \eqref{claim19}, it is convenient to go back to the electric field $\bfe_\eta$ through the second Maxwell equation in \eqref{prb:prin_mu}. As $\eps_\eta=\eps_e$ in 
 $\dom\setminus \Sigma_\eta$, we have:
\begin{eqnarray*} 
\lim_{\eta\to 0} \int_{\dom\setminus \Sigma_\eta}\rt\bsh_\eta\cdot\rt\bphi_\eta &=& \lim_{\eta\to 0} -i\omega\eps_0\int_{\dom\setminus \Sigma_\eta}\bse_\eta(x)\cdot\nabla\rho(x)\wedge\bfv\xeta\\
&=&  -i\omega \eps_0\int_{\dom\times \Sigma^*}\bse_0(x,y)\wedge\nabla\rho(x)\cdot\bfv(y) \ ,
\end{eqnarray*}
 where in the last line we used the weak two-scale convergence $\bse_\eta \cvd \bse_0$.
Recall that by \eqref{PE0}, we have $\bse_0=0$ in $\Sigma$ and $\rt_y\bse_0=0$ in $Y$. 
In particular, for a.e. $x\in\dom$, the periodic vector field $\bff(y) =\bse_0(x,\cdot)\wedge\nabla\rho(x)$ is divergence-free and vanishes in $\Sigma $.
 Therefore by applying \eqref{def:circ}, we obtain $\int_Y \bff(y)\cdot \bfv(y) = (\int_Y \bff)\cdot (\oint\bfv)= 0$.
 The claim \eqref{claim19} follows and the proof of Lemma \ref{variational} is complete. 
\qed

%

\paragraph{Variational characterization of $\bsh_0$.}
For fixed $x\in \dom$, we look for a solution $\bsh_0(x,\cdot)$ to \eqref{weakform}. Let $\bsh(x) =\oint \bsh_0(x,\cdot)$.
 By exploiting the two first equations in \eqref{PH02}, we notice that
the vector field\,  $\bsu:= \bsh_0(x,\cdot) - \bsh(x)$ belongs to 
the following subspace of $\Hunper(Y;\C^3)$
\begin{equation}\label{def_X}
  X_0^{\dv} \ :=\ \Big\{\bfu\in\Hunper(Y;\C^3)\ :\ \rt\bfu=0\ \mbox{in } \Sigma^*\ ,\ \dv\bfu=0 \ \mbox{in } Y\ ,\ \oint\bfu=0 \Big\}\ .
\end{equation}
Accordingly $\bsh_0(x,\cdot)$ has to be searched in $\C^3 \oplus X_0^{\dv}$.
 Notice that the previous sum is direct since non-zero constant functions have a non-vanishing circulation vector
 (see Remark \ref{lem:def_circulation}).
 
 \med
 Now we rewrite the equality \eqref{weakform} which we intentionally restrict to those elements $\bfv\in X_0^{\dv}$ which are divergence-free
 (see Remark \ref{div=0}).
 With $ \bsh_0(x,\cdot)= \bfz + \bfu(y)$, $\bfz= \oint \bsh_0(x,\cdot)$, we find that $\bfu$ solves in $X_0^{\dv}$ the variational problem
 \begin{equation} \label{prb:var}
b_0(\bfu,\bfv)-\eps_r k_0^2\int_Y \bsu\cdot\ov \bsv=\eps_r k_0^2\int_Y \bfz\cdot\ov \bsv\quad,\quad \forall \bsv\in X_0^{\dv}\ .
\end{equation}
where $b_0$ denotes the Hermitian product
\begin{equation} \label{def:b0}
b_0(\bfu,\bsv):=\ \int_\Sigma\rt \bfu\cdot\rt \ov \bsv\quad  \lp = \int_Y\nabla\bfu:\nabla\ov\bfv\rp\ .
\end{equation}
It is easy to check that $X_0^{\dv}$ is a closed subspace of the Hilbert space $\Hunper(Y;\C^3)$ (thanks to \eqref{esti_circ})
and that $b_0$ defined above is a scalar product on $X_0^{\dv}$ which induces an equivalent norm to that $\Hunper(Y;\C^3)$. Indeed we have
\begin{lemma}\label{b0norm}\ 
  There exists a constant $c>0$ such that
  $$   c\, \Vert \bfu\Vert_{\Hunper(Y;\C^3)}^2 \  \le \  b_0(\bfu,\bfu)\ \le\  \Vert \bfu\Vert_{\Hunper(Y;\C^3)}^2 \qquad  
  \text{for all $\bfu\in X_0^{\dv}$}
  $$
  \end{lemma}
\proof
In view of \eqref{equiv_norm_per}, it is enough to show the existence of a constant $k>0$ such that
\begin{equation}\label{Korn}
 \forall \bfv \in X_0^{\dv}\quad,\quad  b_0(\bfv,\bfv)\ =\ \int_Y |\nabla \bfv|^2
\ \ge \ k\,  \int_Y |\bfv|^2\ .
\end{equation}
 By contradiction, assume that \eqref{Korn} does not hold. Then we can find a sequence $(\bfv_n)\subset X_0^{\dv}$ such that
$\Vert \bfv_n\Vert_{L^2}=1$,  $\nabla \bfv_n \to 0$ strongly in $L^2(Y;\C^3)$.
Then, by Rellich's Theorem and possibly after extracting a subsequence,  such a sequence would converge strongly to some constant function  $\bfv$ in $(W_\sharp^{1,2}(Y))^3$ such that $\Vert \bfv\Vert_{L^2}=1$. 
As $X_0^{\dv}$ is a closed subspace of $W_\sharp^{1,2}(Y;\R^3)$,  we need also that this constant function $\bfv$ satisfies $\oint\bfv=0$. This is impossible
unless $\bfv=0$. We get a contradiction with the requirement that $\Vert \bfv\Vert_{L^2}=1$.
 \qed

\begin{lemma}\label{unicite_u}\ Assume that $\Im(\eps_r) >0$.  Then, for every $\bfz\in\C^3$,
 equation \eqref{prb:var}  admits a unique solution in $X_0^{\dv}$\ . \end{lemma}
\proof
\noindent Let $\beta$ a positive real such that $\beta\Im(\eps_r)-\Re(\eps_r)\ge 1$ and let
$$ b(\bfu,\bfv) :=  (1+i \beta)\ \left[b_0(\bfu,\bfv)-\eps_r k_0^2\int_Y \bsu\cdot\ov \bsv\right]\ .$$
Then \eqref{prb:var} is equivalent to solving $ b(\bfu,\bfv) = L_z(v)$ for all $\bfv\in X_0^{\dv}$, being $L_z$ the linear form
on $X_0^{\dv}$ defined by $L_z(v) = (1+i \beta) \eps_r k_0^2 \int_Y \bfz\cdot\ov \bsv$.
Clearly $L_z$ is continuous as well as $b(\cdot,\cdot)$ as a sequilinear form. On the other hand, $b$ is coercive
since, for every $\bfu\in X_0^{\dv}$, it holds
$$
\Re\Big(b(\bsu,\bsu)\Big) =\int_Y\Big(|\nabla\bfu|^2+k_0^2\big(\beta\Im(\eps_r)-\Re(\eps_r)\big)|\bfu|^2\Big) \ge \ b_0(\bsu,\bsu)\ .
$$
The existence and uniqueness of the solution of \eqref{prb:var} follows from Lax-Milgram Lemma.
\qed

\med
Let us notice that the previous existence and uniqueness still holds if $\Im(\eps_r)=0$ provided $\eps_r k_0^2$ does not belong to the discrete set $\{\lambda_n\}$
defined in \eqref{wn}. This a consequence of Fredholm's alternative since the resolvent associated with $b_0$ as an operator on  $L^2_\sharp (Y)$  
 turns out to be compact (see Section \ref{sec:new_form_spec_num}).

\paragraph{Shape magnetic fields.}  By applying Lemma \ref{unicite_u}
 to $\bfz=\bfe_k$ for $k=1,2,3$, we obtain three vector fields $\bfu^1, \bfu^2,\bfu^3$ in $X_0^{\dv}$.
 We associate the following  periodic  fields 
 \begin{equation}\label{decomp_Hk}
\bsh^k(y):=\bfe_k+\bfu^k(y)\quad, \quad \bsj^k := \rt \bsh^k \qquad   k\in\{1,2,3\}\  .
\end{equation}
By linearity,  $\bsh^k$  is  characterized as the unique vector field  $\bfw\in \C^3 \oplus X_0^{\dv}$  satisfying the variational equation \eqref{weakform}
such that $\oint \bfw =\bfe_k$. It is now straightforward to deduce 
 


\begin{proposition}\label{prop:decomp_H0} \ The family $\big\{\bsh^1,\bsh^2,\bsh^3\big\}$ is a basis of solutions 
for equation \eqref{weakform}. Accordingly the two-scale limits  $\bsh_0, \bsj_0$  defined in  \eqref{twoscalconv} 
are uniquely determined in term of  $\disp\bsh(x)$ given in \eqref{effEH} as follows:
\begin{equation}\label{decomp_H0}
\bsh_0(x,y)\ =\ \sum_{k=1}^3 H_k(x)\, \bsh^k(y)\quad , \quad \bsj_0(x,y)\  =\ \sum_{k=1}^3 J_k(x)\, \bsj^k(y) \qquad \text{  a.e. $(x,y)\in \dom\times Y$},
\end{equation} 
being  $ H_k(x) = \bsh(x)\cdot\bfe_k\ $ and\ $J_k:= -\frac1{i\omega \eps_0} H_k $.
\end{proposition}
\proof Let $V$ be the subspace of $\C^3 \oplus X_0^{\dv}$ consisting of all solutions to \eqref{weakform}.  Let $L: V \mapsto \C^3$ the
linear map defined by $L(\bfw):= \oint \bfw$. 
We are reduced to check that this map is bijective. The surjectivity is trivial since, as noticed before, it holds $L(\bsh^k)=\bfe_k$.
On the other hand, if $L(\bfw)=0$, then it means that $\bfw\in X_0^{\dv}$ and that it satisfies \eqref{prb:var} for
 $\bfz=0$. Thanks to the uniqueness result in Lemma \ref{unicite_u}, we conclude that $\bfw=0$.
 \qed
 
 \begin{remark}\label{div=0}\nshape{
If $\bfu\in X_0^{\dv}$  solves \eqref{prb:var}, then the variational equality \eqref{weakform} is satisfied by $\bfw= \bfz +\bfu$ 
 even if the test function $\bfv$ is not divergence-free. Indeed any admissible $\bfv$ for \eqref{weakform}  can be decomposed as $\bfv=\tilde \bfv + \nabla p$ with $\tilde \bfv \in  X_0^{\dv}$ and $p\in \Hunper(Y)$ solving $\Delta p = \dv \bfv$
(as $\int_Y \dv \bfv= 0$ such a $p$ exists and it satsifies  $\oint \nabla p=0$). Therefore, by applying \eqref{prb:var} with $\bfu$ and $\tilde\bfv$, we see that $(\bfw,\tilde \bfv)$ satisfies \eqref{weakform}. In fact this be also the case for $\bfv$ since
$\rt \bfv = \rt \tilde\bfv$ and $\int_Y \bfw\cdot \nabla p =0$ ($\bfw$ is divergence-free like is $\bfu$).
}\end{remark}

\begin{remark}\nshape{
Let us stress that it is necessary that $\Sigma^*$ is \textit{simply-connected} in order to obtain a decomposition of the magnetic field $\bsh_0$  with three independent shape functions like in \eqref{decomp_H0}. In fact the dimension $d$ of the space of solutions to \eqref{prb:var}, which is $d=3$ in our case, can increase with the topological complexity of the inclusion $\Sigma$.  If for instance   $\Sigma$ is a  torus, then $d=4$ as it is shown in \cite{bou_schw} where an extra shape function $\bsh^4$ is needed.
}\end{remark}

\subsection{Integral representation and properties of the effective permeability tensor.} \label{permeabilitytensor}
As will be checked in Section \ref{sec:new_form_spec_num}, the shape functions $\bsh^k$ defined in \eqref{decomp_Hk} coincide with that introduced in \eqref{Hk}.
Therefore, in view of \eqref{bulkmean},  the symmetric effective permittivity tensor $\muef$ given by  \eqref{def_muef_intro} can be recast componentwise through the following integrals
\begin{equation}\label{def_muef}
 \muef_{kl}\ =\ \int_Y\bsh^k\cdot\bfe_l\quad,\quad k,l\in\{1,2,3\}\ .
\end{equation}
The dissipativity property of $\muef$ (needed for the uniqueness issue of Lemma \ref{wellposed})   and   some other  useful relations are
 given in Lemmas \ref{lem:mu} and \ref{quadramean} below.
\begin{lemma}\label{lem:mu}
 Under \eqref{hyp_im_eps_pos}, the real symmetric tensor $\Im(\muef(\omega))$  is positive definite.
 \end{lemma}
 \proof
 Let $\bfz=(z_k)\in \R^3$ and $\bfu^z:=\sum_kz_k\bfu^k$ be the unique solution of  problem \eqref{prb:var}. By using \eqref{decomp_Hk} and \eqref{def_muef}, we get
\begin{equation}\label{mueef_zz}
 \muef\bfz\cdot \bfz=\sum_{k,l=1}^3z_k z_l\int_Y\bsh^k\cdot\bfe_l\ =\ 
\  |\bfz|^2+\int_Y\bfu^z\cdot\bfz .
\end{equation}
 By taking $\bfv=\bfu^z$ in \eqref{prb:var} and dividing by $\eps_r\, k_0^2$,  we infer that
$$
\frac1{\eps_r k_0^2} b_0(\bfu^z,\bfu^z)- \int_Y|\bfu^z|^2\ = \ \int_Y\bfz\cdot\ov\bfu^z\ .
$$
Passing to the conjugate in the equality above and plugging  in \eqref{mueef_zz}, we are led to:
\begin{equation}\label{im_mu_zz}
 \Im(\muef) \bfz \cdot \bfz \  =\ \frac{\Im(\eps_r)}{|\eps_r|^2 k_0^2} \, b_0(\bfu^z,\bfu^z) 
\end{equation}
holding for every $z\in\R^3$. As $b_0$ is a scalar product (see Lemma \ref{b0norm}) and $\Im(\eps_r) >0$,  
we conclude that $\Im(\muef) \bfz \cdot \bfz >0$ unless $\bfu^z=0$ which,
by looking at equation \eqref{prb:var}, is clearly equivalent to $z=0$.
 \qed

\begin{lemma}\label{quadramean} With the notations of Proposition \ref{prop:decomp_H0}, it holds for a.e. $x\in \dom$:
\begin{equation}\label{quadrameanH0}
 \muef \bsh(x) \cdot  \ov\bsh(x)= \int_Y  |\bsh_0(x,\cdot)|^2  - \frac1{\eps_r} \frac{\eps_0}{\mu_0} \int_Y |\bsj_0(x,\cdot)|^2\ ,
\end{equation}
\end{lemma}
\proof We apply Lemma \ref{variational} with $\bfw= \bsh_0(x,\cdot)$ and 
$\bfv = \bar \bsh_0(x,\cdot)- \ov \bsh(x)$ which by construction has a vanishing circulation vector. Then after dividing by $\eps_r\,k_0^2 $ and taking into account
  \eqref{bulkmean},  equation \eqref{weakform} becomes:
$$ \frac1{ \eps_r \,k_0^2} \int_Y |\rt_y \bsh_0(x,\cdot)|^2 \ =\   \int_Y \left( |\bsh_0(x,\cdot)|^2  - \bsh_0(x,\cdot)\cdot \ov\bsh(x)\right) \ =\ 
 \int_Y  |\bsh_0(x,\cdot)|^2 - \muef \bsh(x) \cdot  \ov\bsh(x) \ .$$
Then relation \eqref{quadrameanH0} follows by taking into account the second equation in \eqref{PH0} and relation $k_0^2=\eps_0\,\mu_0\omega^2$.
  \qed

%
%
\section{Spectral description of the effective permeability tensor}\label{sec:new_form_spec_num}

The effective permeability tensor $\muef$ was introduced in Section \ref{sec:reslut_prin_mu} with two different expressions namely \eqref{def_muef_intro2} and \eqref{def_muef_intro}.
The first one is related to  a spectral problem  \eqref{pos_A} set on the Hilbert space $Z_0$  appearing in \eqref{def_Z0} that we recall here:
\begin{equation*}
 Z_0:=\Big\{\bs{f}\in L^2(Y;\R^3)\ :\ \dv \bs{f}=0,\quad \bs{f}=0  \mbox{ in }\Sigma^*\Big\}\ .
\end{equation*}
It turns out that the spectral problem  \eqref{pos_A}  is well suited to numerical approximations (see Section  \ref{subsec:simul_mu}).
The second expression \eqref{def_muef_intro} is related to another spectral problem \eqref{spectralw} where the underlying Hilbert space   $X_0^{\dv}$ 
was constructed in the previous Section (see \eqref{def_X}) 
$$
X_0^{\dv} \ :=\ \Big\{\bfv\in\Hunper(Y;\C^3)\ :\ \rt\bfv=0\ \mbox{in } \Sigma^* \ ,\ \ \dv\bfv=0\ \text{in }Y\ , \ \ \oint\bfv=0 \Big\}\ .
$$
The aim of this construction was to represent the periodic magnetic shape functions $\bsh^k$ which intervene in our main Theorem \ref{t.main}. It turns out that
the entries of tensor $\muef$ can be also recovered from functions $\bsh^k$ by means of relations \eqref{def_muef}.

In this Section we  show that the spectral problems \eqref{pos_A} and \eqref{spectralw} are well-posed and linked together.
Then we  deduce that all definitions given for $\muef$ that is \eqref{def_muef_intro2} , \eqref{def_muef_intro} and \eqref{def_muef} are in agreement.

\subsection{Spectral equivalence.}
Keeping the  notations of section \ref{sec:reslut_prin_mu}, we associate to every element$\bff\in Z_0$ the function $ \Y   \bff$ defined by
\begin{equation}\label{Gammaf}
   \Y   \bff:=\rt\bpsi_f+\frac12\int_\Sigma\bfy\wedge\bff\ ,
\end{equation} 
being $\bpsi_f\in(\Hunper(Y))^3$  the unique solution of \eqref{def:phi_f}, that is: $-\Delta \bpsi_f=\bff\  \mbox{in }Y$ and $\int_Y\bpsi_f=0$.
\begin{lemma}\label{Z0X0} The linear map $\bff \mapsto \Y \bff$  is a bijective from $Z_0$ to $X_0^{\dv}$  and satisfies
$$ \Y (\rt \bsu) = \bsu \quad,\quad \rt(\Y \bff) = \bff \qquad \text{for all $(\bff,\bsu) \in  Z_0\times X_0^{\dv}$} \ .$$
In particular, for all pairs $(\bff,\bfg)\in Z_0^2$ it holds:\quad  $b_0(\Y \bff, \Y \bfg)\, = \, \int_\Sigma \bff \ov \bfg$ .
Therefore $\Y$ is an isometry between Hilbert space $Z_0$ endowed with the scalar product of $L^2(Y)$ and Hilbert space   $X_0^{\dv}$ 
endowed with $b_0$.
\end{lemma}
\proof Let $\bff\in Z_0$ and let us show that $\Y\bff$ belongs to $X_0^{\dv}$. Obviously it is divergence-free.
On the other hand $\rho:= \dv\bpsi_f$ vanishes as the unique periodic solution of $ - \Delta \rho = \dv\bff =0$ with $\int_Y \rho=0$.
It follows that 
\begin{equation}\label{rotGammaf}
\rt (\Y\bff) \, =\, \rt ( \rt( \bpsi_f) \, =\, -\Delta \bpsi_f \, =\, \bff \ .
\end{equation}
form which follows in particular that $\rt \Y\bff= 0$ in $\Sigma^*$. 
It remains to check  the more tricky part, i.e. that $ \oint\Y\bff=0$. 
Let $\bfz\in\C^3$ be arbitrary and  denote by $\bphi_z$ a function in $(\Hunper(Y))^3$ such that $\bphi_z(\bfy)=-\frac12\bfz\wedge\bfy$ in $\Sigma$.
 Then by construction the periodic vector field $\bfg:= \rt \bphi_z + \bfz$ is divergence-free, vanishes in $\Sigma$ and satisfies  $\int_Y\bfg=\bfz$.
 In view of the characterization \eqref{def:circ}, one has therefore
\begin{equation}\label{circ2}
\int_Y\Y\bff\cdot\bfg\ = \ \bfz \cdot \disp \oint\Y\bff \ .
\end{equation}
 By integrating by parts and taking into account that $\ \rt(\rt \bpsi_f)= \bff$, one gets
 $$ \int_Y \rt \bpsi_f \cdot \bfg \,=\, \int_Y \rt \bpsi_f \cdot \rt \bphi_z \,=\, \int_Y \bff \cdot \bphi_z\, =\, - \frac1{2} \, \bfz \cdot \int_\Sigma \bfy\wedge \bff\ .$$ 
Then recalling definition \eqref{Gammaf}, we deduce that $  \int_Y\Y\bff\cdot\bfg =\, 0 \ .$
Thus,  by applying \eqref{circ2}, we are led to $\oint \Y\bff=0$ and  $\Y \bff \in X_0^{\dv}$.

In view of \eqref{rotGammaf}, the Lemma is proved once we have checked that the relation $\Y (\rt \bsu) = \bsu$ holds for every $\bsu\in X_0^{\dv}$.
We observe that, for such $\bsu$, the function $\rt \bsu$ is divergence-free and vanishes in $\Sigma^*$, thus belongs to $Z_0$. Let $\bfv:=\Y (\rt \bsu)$. 
By applying \eqref{rotGammaf} to $\bff= \rt\bsu$, we get the equality $\rt \bsu =\rt \bsv$ from which follows $\bsu=\bsv$ (since $b_0(\bsu-\bsv,\bsu-\bsv)=0$).
\qed 
  
\med
Next we intoduce the sesquilinear form on $Z_0$ defined by
\begin{equation}\label{def:a0}
 a_0(\bff,\bfg)\ :=\  \int_Y \Y \bff \cdot \ov{\Y \bfg} \ =\ \int_Y\nabla\bpsi_f:\ov{\nabla\bpsi_g} \ +\ \frac{1}{4}\ \lp \int_\Sigma \bfy\wedge\bff \rp\cdot
 \lp \int_\Sigma \bfy\wedge \ov \bfg \rp\ ,
\end{equation}
where in the second equality we used definition \eqref{Gammaf}, the fact that the periodic fields $\bpsi_f,\bpsi_g$ are divergence-free (see the proof above) and that their curl have vanishing bulk average over $Y$. 
It is straightforward from Lemma \ref{Z0X0} that $a_0$ is a scalar product on $Z_0$. Moreover it holds 
\begin{equation}\label{continuitya0}
a_0(\bff,\bff) \ \le\ \frac1{k} \, \|\bff\|_{L^2(\Sigma)}^2 \ ,
\end{equation}
where $k$ is the positive constant appearing in \eqref{Korn}.
The spectral problem  \eqref{pos_A} amounts
to finding $(\bff,\alpha)\in Z_0\times\R_+^*$ such that:
 \begin{equation}\label{spec_a0}
 a_0(\bff,\bfg)\,=\, \alpha\, \int_Y \bff\cdot \ov\bfg,\quad\forall \bfg\in Z_0\ .
\end{equation}   
Then,  in relation with the spectral problem \eqref{spectralw}, we obtain the following equivalence principle
holding for every  $(\bsu,\lambda)\in X_0^{\dv}\x\R^+$:
\begin{equation}\label{equiv_spectral}
 b_0(\bsu,\bsv)=\lambda\int_Y \bsu\cdot \ov\bsv,\quad\forall \bsv\in X_0^{\dv}\quad \Longleftrightarrow \quad (\bff,\alpha)= (\rt\bsu, \frac1{\lambda})\quad \text{solves}\quad \eqref{spec_a0}
\end{equation}

\subsection{Reduction to the diagonal form.}\label{ssec:spec}

In view of the equivalence \eqref{equiv_spectral}, we may limit ourselves to the study of the spectral problem \eqref{spec_a0}.
 Since by \eqref{continuitya0} the hermitian product $a_0$ is continuous, we may consider the bounded linear operator $K:Z_0\mapsto Z_0$ such that
 \begin{equation}\label{opA}
 \langle K \bff,\bfg\rangle_{L^2(Y)} \ =a_0(\bff,\bfg)\ ,\qquad\text{for all }(\bff,\bfg)\in Z_0\ .
\end{equation}
We also introduce  $P_0: L^2(Y;\R^3)\mapsto Z_0$ the orthogonal projector on $Z_0$ and  the finite-rank operator $M:Z_0\mapsto L^2(Y;\R^3)$ 
defined for every $\bff\in Z_0$ by setting 
\begin{equation}\label{op:M}
 M\bff(z) \ =\ \frac{1_{\Sigma}(z)}{4}\lp\int_\Sigma\bfy\wedge\bff\rp\wedge\bfz\quad \forall z\in Y\  .
\end{equation}

\begin{lemma}\label{Acompact} With the notations above
 the linear operator  $K: Z_0\to Z_0$ is compact positive and self-adjoint. 
 For all $\bff\in Z_0$, it holds
 \begin{equation}\label{opA2}
K \, \bff\ =\ P_0\lp\bpsi_f+ M \, \bff\rp\ ,
\end{equation}
\end{lemma}


\med
\proof Since $a_0$ is a continuous scalar product,  $K$ is a bounded positive self-adjoint operator on $Z_0$.
Let us take $\bff\in Z_0$ and establish the relation \eqref{opA2}. In view of \eqref{opA}, It is enough to show that for every element $\bfg\in Z_0$, it holds:
\begin{equation}
 a_0(\bff,\bfg) \ =\  \int_Y (\bpsi_f + M \bff)\cdot \ov g \ .\label{proj}
\end{equation}
By \eqref{op:M} one has \ $ \int_Y  M \bff\cdot \ov g  \ = \frac{1}{4}\ \lp \int_\Sigma \bfy\wedge\bff \rp\cdot
 \lp \int_\Sigma \bfy\wedge \ov \bfg \rp\ .$
On the other hand,  by integrating by parts and taking into account that  $\bpsi_g$ satisfies $-\Delta\bpsi_g = \bfg$,  we obtain
$ \int_Y \bpsi_f \cdot \ov g=  \int_Y\nabla\bpsi_f:\ov{\nabla\bpsi_g}$. 
Recalling the expression on the right-hand side appearing in the definition \eqref{def:a0} of $a_0$, we are led to \eqref{proj}.
Thus \eqref{opA2} is established. The compactness property of operator $K$ is then straightforward. Indeed, let  $(\bff_n)$ be a weakly convergent sequence of $Z_0$.
Then $\bpsi_{f_n}$ is bounded $(\Hunper(Y))^3$ hence strongly convergent in $L^2(Y)^3$. It is also the case of $M \bff_n$ since $M$ is finite rank.
The conclusion follows by the continuity of $P_0$ as an operator from $L^2(Y)^3$ to $Z_0$.
\qed

\med
As a consequence of  Lemma \ref{Acompact} and of the equivalent principle \eqref{equiv_spectral}, we have  
 the following result:
\begin{lemma}\label{diagb0}
There exists  a sequence of positive numbers such that $\alpha_0 \ge \alpha_1 \ge \dots \ge \alpha_n$ 
with $\alpha_n\to 0$ and an orthonormal basis  $\{\bff_n\}_{n\in\N}$ of  $Z_0$  consisting of real valued functions such that  
\begin{equation}\label{a0_diag}
 a_0(\bff_n,\bfg)=\alpha_n\int_Y\bff_n\cdot\ov
 \bfg\ ,\qquad\forall n\in\N\quad ,\quad \forall\bfg\in Z_0\ .
\end{equation}
Moreover, by setting for every $n\in \N$
\begin{equation}\label{wn2}
\bsu_n=\frac{1}{\sqrt{\alpha_n}}\ \Gamma\bff_n\ ,\qquad\lambda_n=\frac{1}{\alpha_n}\ ,
\end{equation}
we obtain an orthogonal basis  $\{\bsu_n\}$ of real valued functions in $X_0^{\dv}$ endowed with the scalar product $b_0$ such that for all $n\in\N$:
\begin{equation} \label{b0_diag}
b_0(\bsu_n,\bsv)\ = \, \lambda_n\,  \int_Y  \bsu_n \cdot \ov \bsv  \qquad  \forall \bsv\in X_0^{\dv}\quad,\quad  
 \int_Y  \bfu_n\cdot \bfu_m  =\delta_{nm}.
\end{equation}

\end{lemma}
\proof In view of the compactness property obtained in Lemma \ref{Acompact}, the existence of the pairs $(\bff_n, \alpha_n)$
satisfying \eqref{a0_diag} is obvious. As $\ov K \bff = K \ov \bff$ for every $\bff\in Z_0$, we can chose the eigenvectors $\bff_n$ to be real.
Then we deduce \eqref{b0_diag} by applying the equivalence \eqref{equiv_spectral} to $\bsu_n$ given in \eqref{wn2}.
On the other hand, by the isometry property  in Lemma \ref{Z0X0}, we know that  $\{\Y \bff_n\}$ is orthonormal basis of $(X_0^{\dv},b_0)$. Thus 
$\{\bsu_n\}_{n\in\N}$ is an orthogonal basis of eigenvectors normalized so that 
$b_0(\bsu_n, \bsu_n)\, = \, \lambda_n $ and $\int_Y  |\bfu_n|^2 = 1 \ .$ 
\qed

\subsection{Power series representation of the effective permeability.}\label{subsec:anal_spec_mu}
%

We are now in a position to express  the effective permeability tensor $\muef$ given in  \eqref{def_muef} in terms of a power series expansion.
This is done by projecting  the solution of  the variational equation \eqref{prb:var} on the orthogonal basis $\{\bsu_n\}_{n\in\N}$.
As a byproduct we recover representation formulas \eqref{def_muef_intro2} , \eqref{def_muef_intro} 
making explicit the dependence of $\muef$ with respect to the incident wave number $k_0 = \omega \sqrt{\eps_0 \mu_0}$.

\begin{proposition}\label{prop:reform_mu_num}
Let $(\lambda_n,\bfu_n)$ as in Lemma \ref{diagb0}.
Then the magnetic shape functions $\bsh^k$ introduced in \eqref{decomp_Hk} are given by
\begin{equation}\label{H0_annexe}
 \bsh^k\ =\  \bfe_k + \sum_{n\in\N} \langle\bfe_k,\bfu_n\rangle\frac{\eps_r k_0^2}{\lambda_n-\eps_r k_0^2}\, \bfu_n\quad,\quad  k\in\{1,2,3\}\ .
\end{equation}
Accordingly the  tensor $\muef$ defined in \eqref{def_muef}
admits the following representation:
\begin{equation}\label{mueff_annexe3}
\muef_{ij}(k_0)\ =\
\delta_{ij} + \sum_{n}\frac{\eps_r k_0^2}{\lambda_n-\eps_r k_0^2}\lp \int_Y\bfu_n\cdot \bfe_i\rp\ \lp\int_Y \bfu_n\cdot \bfe_j\rp\ .
\end{equation}
or alternatively in terms of $(\alpha_n,\bff_n)$ by
\begin{equation}\label{mueff_annexe4}
\muef_{ij}(k_0)\ =\ \delta_{ij} + \frac1{4}\, \sum_{n}\frac{\eps_r k_0^2}{1-\eps_r\alpha_n k_0^2}
\lp\int_\Sigma \bs{y}\wedge\bff_n\rp_i\ \lp\int_\Sigma \bs{y}\wedge\bff_n\rp_j \ .
\end{equation}
%

\end{proposition}
\proof  In view of \eqref{def_muef} and \eqref{wn2}, the two series expansions for $\muef$ follow straightforward from  \eqref{H0_annexe}.
Following \eqref{decomp_Hk}, for $k\in \{1,2,3\}$, we write $\bsh^k= \bfe_k + \bsu^k$ where $\bsu^k$ is the unique element of $X_0^{\dv}$ solving
\eqref{prb:var} for $z=\bfe_k$ that is
 \begin{equation} \label{prb:uk}
b_0(\bsu^k,\bsv)-\eps_r k_0^2\int_Y \bsu^k\cdot\ov \bsv \ =\ \eps_r k_0^2\int_Y \bfe_k\cdot\ov \bsv\quad,\quad \forall \bsv\in X_0^{\dv}\ .
\end{equation}
As $\{\bsu_n\}$ is an orthogonal basis of  $X_0^{\dv}$, we may expand $\bsu^k$ as
\begin{equation}\label{decomp_uk}
\bfu^k=\sum_{n\in\N} c_n^k\, \bfu_n\quad,\quad\mbox{where }\quad c_n^k = \int_Y \bsu^k \cdot \bsu_n= \frac1{\lambda_n} b_0(\bsu^k, \bsu_n) \ .
\end{equation}
Note that  the convergence of the series above holds with respect to the Sobolev norm of $\Hunper(Y)$, thus also in $L^2(Y)$.
Plugging  \eqref{decomp_uk} in \eqref{prb:uk} where we take $\bsv= \bsu_n$, we obtain the relations
$$ (\lambda_n-\eps_r k_0^2) \, c_n^k\ =\ \eps_r k_0^2\int_Y \bs{e}_k\cdot\bfu_n\quad ,\quad \ \forall n\in\N \ ,$$
from which follows \eqref{H0_annexe}. \qed

\section{Proof of the main result}\label{sec:proof_prin}

This section is devoted to the proof of Theorem \ref{t.main}. We will proceed according to the following steps:

\begin{itemize}
 \item In a first step, we work under the following energy bound which will be proven \textit{a posteriori} in a second step.
\begin{equation}\label{born_prio_mu}
 \sup_{\eta>0}\int_{B_R}|\bse_\eta|^2+ |\bsh_\eta|^2<+\infty\ ,
\end{equation}
where $B_R$ is an open ball containing $\ov \dom$.

\begin{itemize}
 \item In Proposition \ref{exterieur} we start by establishing the uniform convergence of a subsequence of $(\bse_\eta,\bsh_\eta)$  in every compact subset of $\R^3\setminus\dom$
 from which follows the convergence of the flux of Poynting vectors and a uniform $L^2$ estimate for the rescaled displacement current $\bsj_\eta$
 (see Lemma \ref{flux_pointing}). This estimate together with hypothesis \eqref{born_prio_mu}
 legitimates, up to extracting a subsequence,  the two-scale analysis of the triple $(\bse_\eta,\bsh_\eta,\bsj_\eta)$ performed in Section  \ref{sec:electromagnetic}.

 \item In Proposition \ref{pb-hom2}, we establish  that the \textit{effective} electromagnetic field $(\bse,\bsh)$  deduced from the analysis of Subsections \ref{sec:electric} and \ref{sec:magnetic} and
 of Proposition \ref{exterieur}  satisfies  the limit diffraction problem given in \eqref{prob_lim_mu_intro}. Note that,  thanks to the results of Section \ref{sec:new_form_spec_num}, the tensor $\muef$ defined therein  (see \eqref{def_muef_intro2} or \eqref{def_muef_intro}) agrees with the expression obtained in \eqref{def_muef}. On the other hand the uniqueness of the solution to  \eqref{prob_lim_mu_intro}  (see Lemma \ref{wellposed}) ensures that the whole sequence $(\bse_\eta,\bsh_\eta,\bsj_\eta)$ is weakly two-scale convergent. 
 
 \item In the crucial Subsection \ref{sec_strong}, we improve the convergence above and establish the \textit{strong two-scale} convergence of the triple  $(\bse_\eta,\bsh_\eta, \bsj_\eta)$.
\end{itemize}

\medskip
\item In a second Step, we establish the bound  \eqref{born_prio_mu} by means of a contradiction argument (Lemma \ref{lem:hypo_true}) in which we  use once again the \textit{uniqueness} of the solution of the limit problem  and  the \textit{strong two-scale} convergence  established in the first step.
\end{itemize}

\subsection{Behavior far from the obstacle and energy estimate.}\label{subsec:far}

This first step is related to the convergence of $(\bse_\eta,\bsh_\eta)$  at a positive distance from $\dom$. 
 We use the following result whose proof relies on the hypo-ellipticity of the operator $\Delta + k_0$ (ses for instance \cite{boufel_cap}*{Lemma 2.1}.

\begin{proposition}\label{exterieur}
 Let 
$(\bse_\eta,\bsh_\eta)$ be the solution of problem \eqref{prb:prin_mu} associated to an incident wave $(\bseinc_\eta$, $\bshinc_\eta)$ such that 
 $(\bseinc_\eta,\bshinc_\eta)\to(\bseinc,\bshinc)$ uniformly as $\eta\to 0$. 
If  $(\bse_{\eta},\bsh_{\eta})\cvw(\bse,\bsh)$ weakly in $(L^2(B_R\setminus \dom))^3$, then we can extend $\bse$ and $\bsh$ to $\R^3\setminus \dom$ 
so that both satisfy the Helmholtz equation $\Delta \bsu +k_0^2\bsu=0$ in $\R^3\setminus\ov\dom$. Moreover  the convergence $(\bse_{\eta},\bsh_{\eta})\to(\bse,\bsh)$ holds in $C^\infty(K)$ for all 
compact $K\Subset \R^3\setminus \ov\dom$ while
$(\bse-\bseinc,\bsh-\bshinc)$ satisfies Silver M\"uler condition  \eqref{SM} with respect to $(\bseinc,\bshinc)$. \end{proposition}
A first consequence of the previous result is the  boundedness (see Lemma below) of the flux of the Poynting vector on $\partial B_R$  given by
\begin{equation}\label{def:Pointng} 
 \mP_\eta:=\int_{\partial B_R}\bse_\eta\wedge\ov\bsh_\eta\cdot\bfn\,d\sigma\ .
\end{equation}
It is useful to rewrite  $\mP_\eta$ as a bulk integral. Integrating by parts with the help of the identity 
$$\dv \big(\bse_\eta\wedge\ov\bsh_\eta\big) \ =\ \rt\bse_\eta\cdot\ov\bsh_\eta-\rt\ov\bsh_\eta\cdot\bse_\eta\ $$ 
and taking into account equations \eqref{prb:prin_mu}, we obtain: 
\begin{equation}\label{flux_pointing}
 \mP_\eta=\int_{B_R} (\rt\bse_\eta\cdot\ov\bsh_\eta-\rt\ov\bsh_\eta\cdot\bse_\eta)
\ =\ i\omega\int_{B_R}\left(\mu_0|\bsh_\eta|^2-\eps_0\ov\eps_\eta|\bse_\eta|^2\right )\ .
\end{equation}
The real part of $ \mP_\eta$ is independent of the radius $R$ and represents the energy dissipated by the obstacle. Taking into account the definitions \eqref{def_epseta} \eqref{Jeta}, this dissipation can be recast in terms of the displacement current$\bsj_\eta$ as follows
\begin{equation}\label{energy_pointing}
 \Re\left(\mP_\eta\right)\ =\ -\omega\,\eps_0\,\frac{\im(\eps_r)}{|\eps_r|^2}\int_{\Sigma_\eta} | \bsj_\eta|^2 \ .
\end{equation}

\begin{lemma}\label{bound.jeta}
Let $(\bse_\eta,\bsh_\eta)$ satisfy the uniform upper bound \eqref{born_prio_mu} and solve problem \eqref{prb:prin_mu}. 
Then with $\mP_\eta$  defined by \eqref{def:Pointng} and  $\bsj_\eta$ defined in \eqref{Jeta}, it holds
 \begin{equation}\label{born_Jeta}
 \sup_\eta \left| \mP_\eta\right| < +\infty \quad,\quad \sup_{\eta>0}\|\bsj_\eta\|_{L^2(B_R)}<+\infty\ .
\end{equation}
\end{lemma}
\proof
Up to passing to a subsequence, we can assume that $\limsup_\eta \left|\mP_\eta\right|= \lim_\eta \left|\mP_\eta\right|$
and, by the uniform upper bound \eqref{born_prio_mu}, that $(\bse_{\eta},\bsh_{\eta})$ does converge weakly in $L^2(B_R)$.
We may therefore apply the convergence in  $C^\infty(K)$ obtained in Proposition \ref{exterieur} with $K$ containing a neighbourhood of $\partial B_R$.
This shows that $\limsup_\eta \left|\mP_\eta\right|<+\infty$.
Then under \eqref{hyp_im_eps_pos}, the second upper bound in \eqref{born_Jeta} follows directly from \eqref{energy_pointing}
(notice that in view of \eqref{born_prio_mu} and \eqref{flux_pointing}, we can reach the same conclusion assuming merely that $\eps_r\not=0$).
\qed

%

\subsection{Weak convergence to the solution of the homogenized problem.} \label{weakhom}
Let $\bs{\eps}, \bs{\mu} $ be defined by  \eqref{epsmux}, \eqref{def_epsef_intro}, 
\eqref{def_muef_intro}. Let $(\bse,\bsh)$ be the unique solution of \eqref{prob_lim_mu_intro} and  
let $B_R$ be any ball such that $\dom\Subset B_R$. We associate the fields $\bse_0, \bsh_0$ defined on $B_R\times Y$  by \eqref{E0}\eqref{H0} 
as well as the field $\bsj_0$ determined in Lemma \ref{p.magnetic}. The flux of the Poynting vector through $\partial B_R$ (already used in the proof of Lemma \ref{wellposed})
will be denoted
\begin{equation}\label{effPointing}
 \mP (= \mP(R)) \ :=\  \int_{\partial B_R}\bse\wedge\ov\bsh\cdot\bfn\,d\sigma\ .
\end{equation}

\begin{proposition}\label{pb-hom2}\ \\
Let us assume that the upper bound  \eqref{born_prio_mu} holds and let $(\bse_\eta,\bsh_\eta)$  be the solution  of the diffraction problem \eqref{prb:prin_mu} and
$\bsj_\eta$ given by \eqref{Jeta}.
Then we have the weak two-scale convergence of the triple  $(\bse_\eta,\bsh_\eta,\bsj_\eta) \cvd(\bse_0,\bsh_0,\bsj_0)$ in $B_R$ whereas 
$(\bse_\eta,\bsh_\eta)\to (\bse,\bsh)$ in $C^\infty(K)$ for every compact subset  $K\Subset \R^3\setminus\dom$.
In particular, it holds\, $\mP_\eta\to \mP$ as $\eta\to 0$.


\end{proposition}
\proof
Under the assumption \eqref{born_prio_mu} and thanks to Lemma \ref{bound.jeta} we may, up to passing to a subsequence, assume 
  that the convergences \eqref{twoscalconv} hold for suitable fields $\hat\bse_0$, $\hat \bsh_0, \hat \bsj_0$ belonging to $L^2(B_R\x Y)$ ,  namely
\begin{equation}\label{cvd}
 \bse_{\eta}\cvd \hat\bse_0\ ,\qquad \bsh_{\eta}\cvd \hat\bsh_0\quad\mbox{and}\quad \bsj_{\eta}\cvd \hat\bsj_0\ .
\end{equation}
We are done if we show that the triple $(\hat \bse_0, \hat \bsh_0,\hat \bsj_0)$ agrees with the triple $( \bse_0,  \bsh_0,\bsj_0)$
appearing in our statement (i.e. \eqref{E0} , \eqref{H0} and \eqref{PH0} for $\bsj_0$). 
Let us define:
\begin{equation}\label{averagehat}
\hat\bse(x)=\int_Y\hat\bse_0(x,\cdot)\ ,
\qquad
\hat\bsh(x)=\oint_Y\hat\bsh_0(x,\cdot)\quad\text{for  $x\in\dom$}\ .
\end{equation}
By applying Propositions \ref{prop:decomp_E0} and \ref{prop:decomp_H0}, we infer that 
\begin{equation}\label{e0h0}
 \hat\bse_0(x,y)=\begin{cases}
	      \hat\bse(x)&\hspace{-0.5cm}\text{in }(B_R\setminus \dom)\x Y\\
             \disp \sum_{k=1}^3(\hat\bse(x)\cdot\bfe_k)\bse^k(y)&\text{in }\dom\x Y
             \end{cases}
\ , \quad
 \hat\bsh_0(x,y)=\begin{cases}
	      \hat\bsh(x)&\hspace{-0.5cm}\text{in }(B_R\setminus \dom)\x Y\\
             \disp \sum_{k=1}^3(\hat\bsh(x)\cdot\bfe_k)\bsh^k(y) &\text{in } \dom\x Y
             \end{cases}
\end{equation}
being  $\bse^k,\bsh^k$ the periodic vector fields introduced in \eqref{def:Ei} and \eqref{decomp_Hk} respectively.
Accordingly we have that $(\bse_{\eta},\bsh_{\eta}) \to (\hat\bse,\hat\bsh)$ weakly in $L^2(B_R\setminus\dom)$.
Then, by applying Proposition \ref{exterieur} with a constant incident wave $( \Einc,\Hinc)$, we may extend $(\hat\bse,\hat\bsh)$ to all $\R^3\setminus B_R$ so that it satisfies radiation condition
\eqref{SM} and Maxwell system in the vacuum  outside $\Omega$. Furthermore $(\bse_\eta,\bsh_\eta)\to (\hat\bse,\hat\bsh)$ in $C^\infty(K)$ for every compact subset of $ \R^3\setminus\dom$.
In particular it follows that  $\mP_\eta \to \hat \mP:= \int_{\partial B_R}\hat \bse\wedge\ov{\hat \bsh}\cdot\bfn\,d\sigma.$

\med  Summarizing  we see that all conclusions of Proposition \ref{pb-hom2} hold, provided we  
can establish that $(\hat\bse,\hat\bsh)$ agrees with the solution $(\bse,\bsh)$ to \eqref{prob_lim_mu_intro}.
As by Lemma \ref{wellposed}  this solution is unique, we are thus reduced to check that $(\hat\bse,\hat\bsh)$ satisfies the two first equations in \eqref{prob_lim_mu_intro}.
In fact, as we know already that $(\hat\bse,\hat\bsh)$ satisfies these equations outside $\Omega$,  we may restrict ourselves to $B_R$.

\med
The first equation is readily derived by passing to the distributional limit in the first equation of \eqref{prb:prin_mu}. Indeed the weak limits in $L^2(B_R)$
of $\bse_\eta, \bsh_\eta$  agree with their bulk averages. In view of \eqref{averagehat} and \eqref{def_muef}, there are given respectively by:
$$ \int_Y \hat\bse_0(\cdot,y)\,dy = \hat \bse(x)   \quad,\quad  \int_Y \hat\bsh_0(\cdot,y)\,dy\ =\  \bmu\, \hat\bsh \ .$$
Thus we are led to  
$$
\rt {\hat \bse}\ =\ i\omega\mu_0\, \bmu\, \hat\bsh\quad\mbox{in }  \mathcal{D}'(B_R)\ .
$$
In order to check the second equation, 
we consider $\varphi\in C^\infty_c(B_R)$ and, for fixed $k\in \{1,2,3\}$, the shape function $\bse^k$ defined in \eqref{def:Ei}. Multiplying  the second equation of  \eqref{prb:prin_mu}
by\, $\varphi(x)\bse^k(x/\eta)$ and  integrating by parts the left-hand side  (recall that $\rt_y\bse^k=0$), we get
$$
\int_{B_R} \bs{H}_\eta(x)\cdot\lc\nabla\phi(x)\wedge\bse^k\lp\frac{x}{\eta}\rp\rc=-i\omega\eps_0\int_{B_R} \varphi(x)\bs{E}_\eta(x)\cdot\bs{E}^k\lp\frac{x}{\eta}\rp\ .
$$
By exploiting \eqref{cvd}, the strong two-scale convergence of $\bse^k(x/\eta)$ and the product rule \eqref{product} (and comments below),  we may pass to the limit in the equality above:  
\begin{equation}\label{proof_prin2}
 \int_{B_R\times Y}\hat\bsh_0(x,y) \cdot[\nabla \varphi(x)\wedge\bse^k(y)]\ =\ 
- i\omega\eps_0\eps_e\int_{B_R\times Y}\varphi(x) \hat\bse_0(x,y)\cdot\bse^k(y)\ .
\end{equation}
Next we observe that $\hat\bsh_0(x,\cdot )$ is curl-free in $\Sigma^*$ and $\nabla \varphi(x)\wedge\bse^k$ is divergence-free and vanishes in $\Sigma$.
Therefore, thanks to \eqref{def:circ} and  taking into account \eqref{averagehat}, we are led to 
$$\int_Y \hat\bsh_0(x,y )\cdot[\nabla \varphi(x)\wedge\bse^k(y)]\, dy = \hat \bsh(x) \cdot[\nabla \varphi(x)\wedge\bfe_k]\, .$$ 
On the other hand, in view of definition \eqref{epsmux} and \eqref{epsef2}, we have for a.e. $x\in B_R$:
  $$  \eps_e \, \int_Y \hat \bse_0(x,\cdot)\cdot \bse^k =  \beps(x) \hat \bse(x) \cdot \bfe_k\, .$$
Then by applying Fubini's formula  to the left and right sides of \eqref{proof_prin2},  we are led to
\begin{equation}\label{proof_prin3}
 \int_{B_R}\hat \bsh(x) \cdot[\nabla \varphi(x)\wedge\bfe_k]\ =\ 
- i\omega\eps_0\int_{B_R}\varphi(x) \beps(x) \hat \bse(x) \cdot \bfe_k\ .
\end{equation}
Eventually, thanks to the identity $\dv(e_k\wedge \hat\bsh) = -\rt \hat \bsh \cdot e_k$ and by
integrating by parts the left-hand side of \eqref{proof_prin3}, we infer the following equality
 for all $k\in\{1,2,3\}$ and every test function $\varphi\in C^\infty_c(B_R)$:
$$
\int_{B_R} \rt \hat\bsh\cdot\bfe_k \,\varphi\ =\ -i\omega\eps_0\int_{B_R}\beps \hat \bse \cdot\bfe_k\,\varphi(x) \,
$$
Thus we can conclude that $ \rt{\hat \bsh} = -i\omega\eps_0 \beps \hat \bse$\  holds in  $\mathcal{D}'(B_R)$. 
  The proof of Proposition \ref{pb-hom2} is finished.
\qed

\subsection{Strong two-scale convergence.}\label{sec_strong}
In the following, we keep all assumptions and notations of Subsection \ref{weakhom}. We are going to improve the convergence  result of Proposition \ref{pb-hom2}
by showing that the two-scale convergence of the triple $(\bse_\eta, \bsh_\eta, \bsj_\eta)$ holds strongly. This is done within the three next Lemmas.
We begin with the field $\bsj_\eta$ whose strong convergence is  a straightforward consequence of the fact that $\mP_\eta\to \mP$ (see Proposition \ref{pb-hom2}).
\begin{lemma}\label{cvdf1}
The  following  convergence holds
\begin{equation}\label{cvdf_J0}
\lim_{\eta\to0} \ln\bsj_\eta(x)-\bsj_0\xxeta\rn_{L^2(B_R)}=0 \ .
\end{equation}
\end{lemma}
\proof
As $\Re(\mP_\eta) \to \Re(\mP)$, by \eqref{energy_pointing} and \eqref{pointing_unicite} and recalling that $\bsj_\eta=0$ outside $\Sigma_\eta$ , one has
\begin{multline*}
\lim_\eta 
\int_\dom |\bsj_\eta|^2 \ =\ -\frac{|\eps_r|^2}{\omega\,\eps_0\Im(\eps_r)}\lim_\eta\,  \Re(\mP_\eta)\ =\ -\frac{|\eps_r|^2}{\omega\,\eps_0\Im(\eps_r)} \Re(\mP) \ =
\ \frac{\mu_0}{\eps_0}\frac{|\eps_r|^2}{\Im(\eps_r)}  \int_\dom\Im\Big(\muef\bsh\cdot\ov\bsh\Big)\\
 =\ \int_{\dom\times Y}|\bsj_0|^2 \, dxdy\ , 
\end{multline*}
where in the last equality we used \eqref{quadrameanH0}.
The strong two-scale convergence $\bsj_\eta \cvdf \bsj_0$ follows. On the other hand, in view of the decomposition \eqref{decomp_H0},
 it is easy to check that $\bsj_0(x,x/\eta) \cvdf \bsj_0(x,y)$ ($J_0$ is admissible) and  we are led to \eqref{cvdf_J0}.%
\qed

\begin{proposition}\label{cvdf2} The  following  convergence holds
 \begin{equation}\label{cvdf_h0}
\lim_{\eta\to0}\, \ln \bsh_\eta(x)-\bsh_0\xxeta \rn_{L^2(B_R)}\hspace{-0.5cm}=0\ .
\end{equation}
\end{proposition}
 \proof
Our strategy is based on the argument of compensated compactness in Lemma \ref{lem:divrot}, that we apply on the open set $B_R$ to the pair $(\bsu_\eta,\bsv_\eta)$ defined by
$$  \bsu_\eta \ =\ \bsv_\eta \ :=\ \bsh_\eta(x)-\bsh_0 \xxeta\ .$$
By \eqref{bulkmean},  it is clear that $\bsu_\eta \cvw 0$ weakly in $(L^2(B_R))^3$.
We  need only to prove the following claims:
\begin{equation} \label{Claim1}
a_\eta\ :=\ \dv[\bsh_\eta-\bsh_0(x,x/\eta)]  \to 0  \qquad \text{strongly  in $W^{-1,2}(B_R)$}
\end{equation}
 \begin{equation}\label{Claim2}
\bsb_\eta\ :=\  \rt[\bsh_\eta-\bsh_0(x,x/\eta)]  \to 0  \qquad \text{strongly  in $W^{-1,2}(B_R)$}
\end{equation}
Indeed, under \eqref{Claim1}\eqref{Claim2}, we infer from Lemma \ref{lem:divrot} that  the sequence of non-negative functions $|\bsu_\eta|^2$ converges to zero 
in $\mathcal{D}'(B_R)$. On the other hand, we know from Lemma \ref{exterieur} that  $\bsu_\eta \to 0$ uniformly at a positive distance from $\dom$. By using a trivial localization argument,
we  may readily conclude that $|\bsu_\eta|^2 \to 0$ in $L^1(B_R)$ that is  \eqref{cvdf_h0}.
Notice that, by the bound \eqref{born_prio_mu}, $a_\eta$ and $\bsb_\eta$ are bounded thus weakly convergent to $0$ in $W^{-1,2}(B_R)$. 

\paragraph{Proof of Claim \eqref{Claim1}.} Let $\{\varphi_\eta\}$ be a sequence  such that  $\varphi_\eta\cvw0$ weakly in $W^{1,2}_0(B_R)$. 
In view of Lemma \ref{prop:comp_fort_Hmoins1},
we need only to show the convergence to zero of the duality bracket $\langle a_\eta , \varphi_\eta \rangle$.
Since $\dv\bsh_\eta=0$, we have 
\begin{equation}\label{aeta}
\langle a_\eta , \varphi_\eta\rangle =-\int_{B_R}\bsh_0\xxeta\cdot\nabla\varphi_\eta(x)\dx\ .
\end{equation}
Up to passing to a subsequence, we can assume that $\nabla \varphi_\eta \cvd \xi_0$ for a suitable element of $(L^2(\dom \times Y))^3$.
Furthermore, by Proposition   \ref{example_cvdf}  (see \eqref{conv_phi_eta_tilde_bis}), there exists 
$\psi_0\in L^2(B_R; W^{1,2}_\sharp(Y))$ such that $\xi_0 =\nabla_y \psi_0\, .$ 
As  $\bsh^0 = \sum_k H_k(x) \bsh^k(y)$ is a finite sum of admissible functions (see the comments after \eqref{product}), we may pass to the limit in \eqref{aeta}
applying \eqref{product}) for $\varphi=1$. We obtain: 
$$
\lim_{\eta\to 0}  \langle a_\eta , \varphi_\eta \rangle =  - \int_{B_R} \left(\int_{ Y} \bsh_0(x,\cdot)\cdot \nabla_y \psi_0 (x,\cdot)\,dy\right) \, dy\ .
$$
Recalling that  $\bsh_0(x,\cdot)$ is divergence-free (see \eqref{PH0}), we see that the right-hand side vanishes by integrating by parts over $Y$
and by taking into account the periodicity of $\bsh_0(x,\cdot)$ and $\psi_0(x,\cdot)$.

\bigskip\noindent
\paragraph{Proof of Claim \eqref{Claim2}.}
This is the most  delicate part. In a similar way as before, we consider a sequence $\{\bphi_\eta\}$  such that  $\bphi_\eta\cvw0$ weakly in $(W^{1,2}_0(B_R))^3$.
We need  to show that $\langle \bsb_\eta , \bphi_\eta\rangle \to 0$ where
\begin{eqnarray}\label{conv_forte_temp}
\langle \bsb_\eta , \bphi_\eta\rangle &=&\int_{B_R}\rt\bsh_\eta\cdot\bphi_\eta-\int_{B_R}\bsh_0\xxeta\cdot\rt\bphi_\eta  \nonumber\\
&=& -i\omega\eps_0\frac{1}{\eta}\int_{B_R}\bsj_\eta\cdot\bphi_\eta-\int_{B_R}\bsh_0\xxeta\cdot\rt\bphi_\eta\ .
\end{eqnarray}
Here, in the second equality, we used equation \eqref{prb:prin_mu} and the definition \eqref{Jeta}.
It is now  useful to consider the piecewise constant function 
$$[\bphi_\eta]_\eta(x):=\frac{1}{\eta^3}\sum_{k\in I_\eta}\lp\int_{Y^k_\eta}\bphi_\eta\rp1_{Y^k_\eta}(x)\quad,\quad Y_\eta^k=\eta(Y+\bfk)\ ,\ I_\eta:=\{\bfk\in\Z^3\ :\ Y_\eta^k\subset B_R\}\ .$$
By applying Proposition \ref{example_cvdf} to each component of $\bphi_\eta$, we are ensured of the existence
of an element $\bpsi_0\in W^{1,2}(B_R;\Hunper( Y;\C^3))$ such that
$$
\frac{\bphi_\eta-[\bphi_\eta]_\eta}{\eta}\cvd \bpsi_0\quad,\quad \rt_y\bphi_\eta\cvd\rt_y\bpsi_0\ .
$$
In view of \eqref{conv_forte_temp}, let us write \
$\langle \bsb_\eta , \bphi_\eta\rangle= r_\eta  +s_\eta$\ where 
\begin{eqnarray}\label{reta}
r_\eta \ &:=&\  - i\omega\eps_0\int_{B_R}\bsj_\eta\cdot\frac{\bphi_\eta-[\bphi_\eta]_\eta}{\eta}- \int_{B_R}\bsh_0\xxeta\cdot\rt\bphi_\eta\ ,\\
s_\eta \ &:=&\   -i\omega\eps_0\int_{B_R}\frac{\bsj_\eta}{\eta}\cdot [\bphi_\eta]_\eta\ . \label{seta}
\end{eqnarray}
Thanks to the strong two-scale convergences $\bsh_0\xxeta \cvdf \bsh_0$ and $\bsj_\eta \cvdf \bsj_0$ (see Lemma \ref{cvdf1}),
 we may apply twice the multiplication rule \eqref{product} (and subsequent comments). We are led to
$$
\lim_{\eta\to 0} r_\eta \ =\ - \int_{B_R\x Y} \left[i\omega\eps_0\bsj_0\cdot\bpsi_0 + \bsh_0\cdot\rt_y\bpsi_0 \right] \, \dx\dy\ .
$$
Here we used the fact that $\bsh^0$ is an admissible function  and also that $\bsj_\eta$ is compactly supported in $B_R$ (so that \eqref{product} can be applied with a localizing test function $\varphi$
such that $\varphi=1$ in $\dom$).
Eventually, as   $\rt_y\bsh_0(x,\cdot) =- i\omega\eps_0\, \bsj_0(x,\cdot)$ (see \eqref{PH0}), by integrating by parts on $Y$
and taking into account the periodicity of $\bpsi_0(x,\cdot)$, we readily deduce that the right-hand member in the equality above vanishes, thus $r_\eta \to 0$.

\med
In order to conclude, it remains to show that $s_\eta$ defined in \eqref{seta} vanishes as $\eta\to 0$. To that aim we consider, 
for $k\in\{1,2,3\}$, a function $\theta_k\in C^\infty_c(Y)$ verifying $\nabla\theta_k=\bfe_k$ in $\Sigma$ and extended by periodicity in $\R^3$.
Then we consider the potential $w_k$ defined by
$$ w_\eta (x)\ :=\ \eta \, \sum_{k=1}^{k=3} \left([\bphi_\eta]_\eta(x) \cdot\bfe_k\right)\ \theta_k\xeta\ .$$
By construction and recalling that $[\bphi_\eta]_\eta$ is piecewise constant, it holds
\begin{equation}\label{Sigmaweta}
\nabla w_\eta \ =\ [\bphi_\eta]_\eta \qquad \text{on $\Sigma_\eta$}\ .
\end{equation}
On the other hand, since the $\theta_k$'s are Lipschitzian , one checks easily  that $w_\eta$ belongs to $W^{1,2}_0(B_R)$ and satisfies $|\nabla w_\eta| \le C \, \left|[\bphi_\eta]_\eta\right|$ for a suitable constant $C>0$. In particular as $[\bphi_\eta]_\eta \to 0$ strongly in $L^2(B_R)$, we have that
\begin{equation}\label{estiweta}
 \lim_{\eta\to 0}  \| \nabla w_\eta - [\bphi_\eta]_\eta \|_{L^2(B_R)} \ =\ 0\ .
\end{equation}
Thanks to \eqref{Sigmaweta} and noticing that $\frac{\bsj_\eta}{\eta}$  agrees with $\eps_e\, E_\eta$ in $B_R\setminus\Sigma_\eta$, we have
$$ s_\eta\ =\ \frac1{\eta} \int_{B_R} \bsj_\eta \cdot \nabla w_\eta + \eps_e  \int_{B_R} \bse_\eta \cdot \left([\bphi_\eta]_\eta - \nabla w_\eta\right) \ .$$
Clearly the first integral above vanishes (since $\bsj_\eta$ is divergence-free) while the second one converges to zero by Cauchy-Schwarz inequality  and by taking into account 
 \eqref{estiweta} and the upper bound \eqref{born_prio_mu}.
Summarizing, we have proved that $\langle \bsb_\eta , \bphi_\eta\rangle= r_\eta +s_\eta \to 0$ and our claim \eqref{Claim2} follows.
The proof of Proposition \ref{cvdf2} is finished.
\qed

\begin{proposition}\label{cvdf3} 
The  following  convergence holds
 \begin{equation}\label{cvdf_e0}
\lim_{\eta\to0} \ln \bse_\eta(x)-\bse_0\xxeta \rn_{L^2(B_R)}\hspace{-0.5cm}=0\ .
\end{equation}
\end{proposition}
\proof 
First we observe that  $|\bse_\eta|^2 = \frac{\eta^2 }{|\eps_r|^2}|\bsj_\eta|^2$ holds in $\Sigma_\eta$ and then, thanks to the strong two-scale convergence of $\bsj_\eta$ obtained in Lemma \ref{cvdf1},  we have
\begin{equation}\label{strongCVE}
\lim_{\eta\to 0} \frac1{\eta^2} \, \int_{\Sigma_\eta}|\bse_\eta|^2\ =\ \frac{1 }{|\eps_r|^2} \int_{\dom\times Y} |\bsj_0|^2 \ .
\end{equation}
 In particular, it holds $ \int_{\Sigma_\eta}|\bse_\eta|^2\to0$. Thus , recalling that $\bse_0(x,\cdot)= \bse(x)$ for $x\in B_R\setminus\dom$ and that $\bse_0 =0$ in
 $\dom\times \Sigma$,
 the convergence \eqref{cvdf_e0} is proved once we have shown that
\begin{equation} \label{limsup1}
\limsup_\eta \int_{B_R\setminus\dom} |\bse_\eta|^2 \ \le\ \int_{B_R\setminus\dom} |\bse|^2 \quad,\quad 
 \limsup_\eta \int_{\dom\setminus\Sigma_\eta} |\bse_\eta|^2 \ \le\ \int_{\dom\times Y} |\bse_0|^2 \ .
\end{equation}
As in Lemma \ref{cvdf1}, the main argument is the convergence $\mP_\eta\ \to \mP$ established in Proposition \ref{pb-hom2}.
It is convenient to rewrite Claim \eqref{limsup1} in the following equivalent form
 \begin{equation} \label{limsup2}
\alpha := \limsup_\eta \left(\int_{B_R\setminus\dom} |\bse_\eta|^2 + \eps_e \int_{\dom\setminus\Sigma_\eta} |\bse_\eta|^2\right) \ \le\  \int_{B_R\setminus\dom} |\bse|^2  
+ \eps_e \int_{\dom\times Y} |\bse_0|^2 \ =\ \int_{B_R} \beps \, \bse\cdot \ov \bse\ ,
\end{equation}
 where in the last equality we used \eqref{epsmux} and \eqref{epsef2}. It is trivial that \eqref{limsup1} implies \eqref{limsup2}. The converse
 implication follows from the fact that one has
 $$\liminf_\eta \int_{B_R\setminus\dom} |\bse_\eta|^2 \ge \int_{B_R\setminus\dom} |\bse|^2\quad,\quad 
  \liminf_\eta \int_{\dom\setminus\Sigma_\eta} |\bse_\eta|^2 \ge \int_{\dom\times Y} |\bse_0|^2.$$

\med 
Now, by taking the imaginary parts in \eqref{flux_pointing} and after dividing by $\omega\eps_0$, we obtain the relation
$$ \int_{B_R\setminus\dom} |\bse_\eta|^2 + \eps_e \int_{\dom\setminus\Sigma_\eta} |\bse_\eta|^2 \ =\ - \frac{\Im( \mP_\eta)} {\omega \eps_0}
 + \frac{\mu_0}{\eps_0} \int_{B_R} |\bsh_\eta|^2 -  \frac{\Re( \eps_r)}{\eta^2}  \int_{\Sigma_\eta}|\bse_\eta|^2
\ .
$$
Letting $\eta\to 0$, with the convergence of $\mP_\eta$ to $\mP$ and thanks to \eqref{strongCVE} and the strong two-scale convergence $\bsh_\eta \cvdf \bsh$,
 we see that the left-hand member $\alpha$ of \eqref{limsup2} is actually a limit and it holds
 \begin{eqnarray}
\alpha\ &=& - \frac{\Im( \mP)} {\omega \eps_0}
 + \frac{\mu_0}{\eps_0} \int_{B_R\times Y} |\bsh_0|^2 -  \frac{\Re( \eps_r)}{|\eps_r|^2} \int_{\dom\times Y} |\bsj_0|^2  \nonumber\\
 \label{final1} &=&\ - \frac{\Im( \mP)} {\omega \eps_0}
 + \frac{\mu_0}{\eps_0} \int_{B_R}\Re\Big({\bs{\mu}}(x)\bsh\cdot\ov\bsh\Big) \ , 
\end{eqnarray}
where $\bmu$ is define in \eqref{epsmux} and in the last equality we took  the integral over $\dom$ of the real parts of the terms appearing in the equality \eqref{quadrameanH0} while for $x\in B_R\setminus\dom$ 
we used the fact that $\bsh_0(x,\cdot) =\bsh(x)$ and $\bsj_0(x,\cdot)=0$.
Now it remains  to evaluate $\Im( \mP)$. In view of the second equality in \eqref{balance} in which we take the imaginary parts and where $\beps$ is real, we get
\begin{equation}\label{final2}
\Im(\mP)=\omega\lp \mu_0\int_{B_R}\Re\Big({\bs{\mu}}(x)\bsh\cdot\ov\bsh\Big)-\eps_0\int_{B_R}\bs{\eps}(x)\bse\cdot\ov\bse \rp\ . 
\end{equation}
Putting \eqref{final1} and \eqref{final2} together, we can conclude that $\alpha = \int_{B_R} \beps \, \bse\cdot \ov \bse$.
Thus the claim \eqref{limsup1} holds and the proof of Proposition \ref{cvdf3} is finished.
%
\qed

\subsection{End of the proof.} 

Collecting the results of Proposition \ref{pb-hom2}, Lemma \ref{cvdf1}, Proposition \ref{cvdf2} and Proposition \ref{cvdf3}, we see that the proof of our main Theorem \ref{t.main} is achieved
once the energy bound \eqref{born_prio_mu} has been established.
The following final Lemma allows us to conclude.
\begin{lemma}\label{lem:hypo_true}
Under the assumptions of Theorem \ref{t.main},  the energy bound given in \eqref{born_prio_mu} holds.
\end{lemma}
\proof
Let us assume by contradiction that 
 $(\bse_\eta,\bsh_\eta)$ solving \eqref{prb:prin_mu} is  not uniformly bounded in $L^2(B_R)^3$. Then it exists
 a subsequence (still denoted $\eta$) such that $\lim_\eta (\int_{B_R}|\bse_\eta|^2+ \int_{B_R}|\bsh_\eta|^2 )=+\infty$. For such a subsequence we denote
\begin{equation}\label{normalis_mu}
 t_\eta:=\lp\int_{B_R}|\bse_\eta|^2+ \int_{B_R}|\bsh_\eta|^2   \rp^{\frac{1}{2}}\to +\infty\quad,\quad \hat\bse_\eta:=\frac{\bse_\eta}{t_\eta}\quad\mbox{and}\quad
\hat\bsh_\eta:=\frac{\bsh_\eta}{t_\eta}\ .
\end{equation}
By construction $(\hat\bse_\eta,\hat\bsh_\eta)$ satisfies \eqref{born_prio_mu} whereas, by linearity, it is the solution of \eqref{prb:prin_mu} when the incident wave 
is given by $(\frac{\bseinc}{t_\eta},\frac{\bshinc}{t_\eta})$. We can then apply Proposition \ref{pb-hom2}: the limit of $(\hat\bse_\eta,\hat\bsh_\eta)$ 
is characterized in term of  the unique solution $(\hat\bse,\hat\bsh)$ to the diffraction problem \eqref{prob_lim_mu_intro}
in which the incident wave $(\hat \bseinc, \hat \bshinc)= \lim_{\eta\to 0} (\frac{\bseinc}{t_\eta},\frac{\bshinc}{t_\eta})= (0,0).$
Therefore, by the uniqueness Lemma \ref{wellposed}, it holds  $(\hat\bse,\hat\bsh)=(0,0)$ and the two-scale limit $(\hat\bse_0,\hat\bsh_0)$ of $(\hat\bse_\eta,\hat\bsh_\eta)$ vanishes as well on $B_\R\times Y$.
Then, by applying the strong convergence results \eqref{cvdf_h0} and \eqref{cvdf_e0}, we are led to 
$$\lim_{\eta\to0}\int_{B_R}|\hat\bse_\eta|^2\ =\ \lim_{\eta\to0}\int_{B_R}|\hat\bsh_\eta|^2\, =\, 0\ .$$
This is impossible  since,  by \eqref{normalis_mu}, we have \ 
$\int_{B_R}|\hat\bse_\eta|^2+\int_{B_R}|\hat\bsh_\eta|^2\ =1\ .$
\qed


\addcontentsline{toc}{section}{References}

\bibliography{biblio}
\bibliographystyle{plain}

\end{document}